\documentclass{article}
\usepackage{amsmath}
\usepackage{amsfonts}
\usepackage[dvips]{graphics}

\newcommand{\NI}{\noindent}

\newcommand{\eps}{\epsilon}

\newcommand{\sn}{\sqrt{n}}
\newcommand{\Th}{\theta}

\newcommand{\tausq}{{\tau}^{2}}
\newcommand{\tautsq}{{\tau_{t}}^{2}}
\newcommand{\taur}{\tau_t^{-2}}
\newcommand{\IN}{\frac{1}{n}}

\newcommand{\be}{\begin{equation}}
\newcommand{\ee}{\end{equation}}
\newcommand{\beq}{\begin{eqnarray}}
\newcommand{\beqn}{\begin{eqnarray*}}
\newcommand{\eeq}{\end{eqnarray}}
\newcommand{\eeqn}{\end{eqnarray*}}
\newcommand{\bed}{\begin{description}}
\newcommand{\ed}{\end{description}}

\newcommand{\ra}{\rightarrow}
\newcommand{\lra}{\longrightarrow}
\newcommand{\cp}{\stackrel{P}{\longrightarrow}}
\newcommand{\cd}{\stackrel{d}{\longrightarrow}}

\newcommand{\tb}{\textbf}
\newcommand{\mb}{\mathbf}

\newcommand{\rt}{\right}
\newcommand{\lt}{\left}

\begin{document}
\begin{center}
\large{Project Report on Resampling in Time Series Models}

\medskip
Submitted by: Abhishek Bhattacharya

\smallskip
Project Supervisor: Arup Bose\\
\end{center}

\noindent
\begin{abstract}
This project revolves around studying estimators for parameters in different Time Series models and
studying their assymptotic properties. We introduce various bootstrap techniques for the estimators
obtained. Our special emphasis is on Weighted Bootstrap. We establish the consistency of this scheme
in a AR model and its variations. Numerical calculations lend further support to our consistency
results. Next we analyze ARCH models, and study various estimators used for different error 
distributions. We also present resampling techniques for estimating the distribution of the 
estimators. Finally by simulating data, we analyze the numerical properties of the estimators.
\end{abstract}

\section{Bootstrap in AR(1) model} 

Let ${X_t}$ be a stationary AR(1) process, that is,
\begin{eqnarray} 
 X_t = \theta X_{t-1} + Z_t \hspace{0.1in} for \hspace{0.1in}t = 1,2,\ldots \label{ar} \\
 Z_t \hspace{0.1in}iid\hspace{0.1in}(0,\sigma^2);\hspace{0.1in}EZ^4_t <\infty ; \hspace{0.1in}|\theta|<1. \nonumber
 \end{eqnarray}

\noindent 
We have assumed $\sigma$ to be known, and $\theta$ is the unknown parameter of interest. Then the 
Least Squares estimate for $\theta$ (which is approximately the MLE in case of normal errors) is given by 
\[
\hat{\theta}_n = \frac{\sum_{t=2}^n X_t X_{t-1}}{\sum_{t=2}^n X_{t-1}^2}
\] 
Then it can be established that 
\begin{equation} \label{ars}
\sqrt{n}(\hat{\theta}_n - \theta)
\stackrel{d}{\longrightarrow}N( 0, (1-{\theta^2}))
\end{equation}
 Let us introduce two particular bootstrap techniques specially used to estimate the distribution of $\hat\theta_n$ from a realization of model (\ref{ar}).\\ \\
\textbf{(a) Residual Bootstrap} Let $\tilde{Z}_t = X_t -  \hat{\theta}_n X_{t-1} ,\hspace{0.1in} t=2,3,\ldots,n$ 
 and let $\hat{Z_t}$ be the standardized version of $\tilde{Z}_t$ such that 
$\frac{1}{n-1}\sum \hat Z_t\;=\;0$ and $\frac{1}{n-1}\sum \hat Z^2_t\;=\;1$.
 Now we draw $Z^*_t$, $t=1,2,\ldots,N$ with replacement from $\hat{Z_t}$ and define
\begin{eqnarray*}
X^*_1 &=& Z^*_1 \\
X^*_t &=& \hat{\theta}_n X^*_{t-1}+Z^*_t\hspace{0.1in}, t=2,\ldots,N.
\end{eqnarray*}
 and form the statistic 
 \begin{equation} \label{rb}
 \hat\theta^*_n = \frac{\sum_{t=2}^n X^*_t X^*_{t-1}}{\sum_{t=2}^n (X^*_{t-1})^2}
\end{equation}
 Then (\ref{rb}) forms an estimator of  $\hat{\theta}_n$ and is called the \textbf{Residual Bootstrap estimator}.
 We repeat the simulation process several times to estimate the distribution of $\hat{\theta}^*_n$.\\ \\ 
\textbf{(b) Weighted Bootstrap} Alternatively we define our resampling estimator 
\begin{equation} \label{wb1}
 \hat{\theta}^*_n= \frac{\sum_{t=2}^n w_{nt} X_t X_{t-1}}{\sum_{t=2}^n w_{nt}(X_{t-1})^2}
  \end{equation}
   where $\{ w_{nt} ; 1 \le t \le n , n \ge 1 \}$ is a triangular sequence of random variables,
independent of $\{X_t\}$. These are the so called ``Bootstrap weights'', and the 
estimator (\ref{wb1}) is the \textbf{Weighted Bootstrap Estimator}.

\subsection{A Bootstrap Central limit theorem}
Under suitable conditions on the weights to be stated below, we establish the distributional consistency of the Weighted Bootstrap Estimator,
\hspace{.05in}$\hat{\theta}^*_n$ defined in (\ref{wb1}). To establish consistency, we will prove a Bootstrap CLT for which we will
need the following established results: 
\newtheorem{res}{Result}
\begin{res}[P-W theorem; see Praestgaard and Wellner(1993)] \label{pw} 
Let \linebreak $\{c_{nj};\; j=1,2,\ldots,n;\; n \ge 1\}$ be a triangular array of constants,
 and let $\{U_{nj} \hspace{0.1in} j=1,2,\ldots,n;\; n \ge 1\}$ be a triangular array of row
   exchangeable random variables such that as $n \rightarrow \infty$,
   \begin{enumerate} 
\item $\frac{1}{n} \sum_{j=1}^n c_{nj} \rightarrow 0$  
\item $ \frac{1}{n}\sum_{j=1}^n c^2_{nj} \rightarrow \tau^2$
\item $\frac{1}{n} {\max}_{1\le j \le n} c^2_{nj} \rightarrow 0 $
\item  $E(U_{nj}) = 0  \hspace{.1in}     j=1,2,\ldots ,n.\hspace{.1in}n\ge 1$ 
\item $E(U^2_{nj}) = 1  \hspace{.1in}     j=1,2,\ldots ,n.\hspace{.1in}n\ge 1$
\item $\frac{1}{n} \sum_{j=1}^n U^2_{nj}\stackrel{P}{\rightarrow}1$
\item ${\lim}_{k\rightarrow \infty}  {\limsup}_{n\rightarrow \infty} \sqrt{(E(U^2_{nj}I_{\{|U_{nj}|>k\}})}= 0 $
\end{enumerate}
Then under the above conditions,
\begin{equation}\label{pwr1}
\frac{1}{\sqrt{n}} \sum_{j=1}^n c_{nj}U_{nj} \stackrel{d}{\longrightarrow}N(0, \tau^2)
\end{equation}
\end{res}

\noindent
Result (\ref{pw}) can be generalized by taking \{$c_{nj}$\} random variables, independent of 
$\{U_{nj}\}$ and the conditions (1), (2) and (3) replaced by convergence in probability. In that case conclusion (\ref{pwr1}) is replaced by
\begin{equation}\label{pwr2}
P \left [\frac{1}{\sqrt{n}} \sum_{j=1}^n c_{nj}U_{nj} \in C \; \vline \{c_{nj}; j=1,\ldots,n; n \ge 1\} \right ] - 
P\left [Y \in C \right ] = o_P(1) 
\end{equation}
where $Y \sim N(0, \tau^2)$ and $C \in \mathcal{B}(\mathbb{R})$ such that $P(Y \in \partial C) = 0$.

\begin{res} \label{res2}
 Let $\{X_1, X_2,\ldots,X_n\}$ be the realization of the stationary AR(1) process (\ref{ar}). Then \\ 
$\frac{1}{n} \sum_{t=1}^{n-k} X^a_t Z^b_{t+k} \stackrel{a.s.}{\rightarrow}E(X^a_t Z^b_{t+k})$ whenever $EZ^{\max{(a,b)}}_t$ $< \infty \;\;$ $\forall a,b,k \in \cal{Z^+}$; $a,b \ge 0\;;\; k>0$. This can be established using the Martingale SLLN; see Hall and Heyde 1980.
 \end{res}
Let us use the notations $P_B ,\; E_B ,\; V_B$ to respectively denote probabilities, expectations and
 variances with respect to the distribution of the weights, conditioned on the given data $\{X_1,\ldots,X_n\}$. The weights are
 assumed to be row exchangeable. We henceforth drop the first suffix in the weights $w_{ni}$ and denote it by $w_i$. 
  Let $\sigma^2_n = V_B(w_i),\;\; W_i = \sigma_n^{-1}(w_i - 1)$. The following conditions on the row
   exchangeable weights are assumed: 
   \begin{description}
\item[A1.] $E_B(w_1) = 1$ 
\item[A2.] $0 < k < \sigma^2_n = o(n)$ 
\item[A3.] $c_{1n}  = Cov(w_1,w_2) = O(n^{-1})$ 
\item[A4.] Conditions of Result(\ref{pw}) hold with $U_{nj}=W_{nj}$.
\end{description}
 \newtheorem{thm}{Theorem}
\begin{thm} \label{t1}
 Under the conditions (A1)-(A4) on the weights,
 \begin{equation}
P_{B}\left [\sqrt{n}\sigma^{-1}_n(\hat\theta^{*}_n - \hat\theta_n) \le x |X_1,\ldots,X_n\right ] - P\left [Y \le x \right ] = o_P(1)\hspace{.1in} \forall x \in \mathbb{R}
\end{equation}
where $Y \sim N(0, (1-\theta^2))$.
\end{thm}

\smallskip 
\noindent 
\textbf{Proof}
  Note that  
\begin{eqnarray*}
\hat\theta^*_n &=& \frac{\sum_{t=2}^n w_t X_t X_{t-1}}{\sum_{t=2}^n w_{t}X^2_{t-1}} \\
&=& \frac{\sum_{t=2}^n w_t X_{t-1}(\theta X_{t-1}+ Z_t)}{\sum_{t=2}^{n} w_t X^2_{t-1}} \\
&=& \theta + \frac{\sum w_t X_{t-1}Z_t}{\sum w_t X^2_{t-1}}
\end{eqnarray*}
Similarly 
\[
  \hat{\theta}_n =\frac{\sum X_t X_{t-1}}{\sum X^2_{t-1}} = \theta + \frac{\sum  X_{t-1}Z_t}{\sum  X_{t-1}} 
\]
Hence 
\begin{eqnarray*}
\hat\theta^*_n-\hat{\theta}_n &=& \frac{\sum w_t X_{t-1} Z_t}{\sum w_t X^2_{t-1}} - 
\frac{\sum  X_{t-1} Z_t}{\sum  X^2_{t-1}} \\
&=& \frac{\sum w_t X_{t-1} Z_t}{\sum w_t X^2_{t-1}} - \frac{\sum  X_{t-1} Z_t}{\sum w_t X^2_{t-1}}
     +\frac{\sum X_{t-1} Z_t}{\sum w_t X^2_{t-1}} - \frac{\sum  X_{t-1} Z_t}{\sum  X^2_{t-1}} \\
&=& \frac{\sum (w_t-1) X_{t-1} Z_t}{\sum w_t X^2_{t-1}}-\frac{\sum X_{t-1}Z_t \sum (w_t-1) X^2_{t-1}}{\sum  X^2_{t-1} \sum w_tX^2_{t-1}}    
\end{eqnarray*}
Now using Result (\ref{res2}),
\begin{eqnarray}
\frac{\sum  X_{t-1} Z_t}{n} &\stackrel{a.s.}{\longrightarrow}& E(X_{t-1} Z_t)\;=\; 0 \label{a1} \\
\frac{\sum  X^2_{t-1} Z^2_t}{n}&\stackrel{a.s.}{\longrightarrow}& E(X^2_{t-1} Z^2_t) \nonumber \\
&=& \sigma^4(1-\theta^2)^{-1} \label{a2}
\end{eqnarray}
\noindent
\textbf{Claim 1.} For $\tau^2 \;= \;\sigma^4(1-\theta^2)^{-1}$,
\[
P_B\left[\;\frac{1}{\sqrt{n}}\sum_{t=2}^n W_t X_{t-1}Z_t \le x\;\vline \;X_1,\ldots, X_n\;\right] \stackrel{P}{\longrightarrow}\Phi(\;\;\frac{x}{\tau}\;\;) \hspace{.1in}\forall x \in \mathbb{R}
\]
To see this let us verify the conditions of Result(\ref{pw}) with $c_{nj} = X_jZ_{j+1}$ and $U_{nj} = W_j$ for $j = 1,\ldots,{n-1}$.
\begin{enumerate} 
\item $\frac{1}{n} \sum_{t=2}^n X_{t-1}Z_t\stackrel{P}{\longrightarrow}0$ 

\smallskip
Follows from (\ref{a1}).

\item $ \frac{1}{n}\sum_{t=2}^n X^2_{t-1}Z^2_t\stackrel{P}{\rightarrow}\sigma^4(1-\theta^2)^{-1} (=\tau^2)$

\smallskip
Follows from (\ref{a2}).

\item $n^{-1}\max(X^2_{t-1}Z^2_t)\stackrel{P}{\longrightarrow}0$
\end{enumerate}
\textbf{Proof} Let $Y_t = X^2_{t-1}Z^2_t = X^2_{t}X^2_{t-1}-2 \theta X_t$
$X^3_{t-1}+\theta^2X^4_{t-1}$

\smallskip
\noindent
Then given $\epsilon > 0$,
\begin{eqnarray*}
& & P(n^{-1}\max Y_t >\epsilon)\;=\; P(\max Y_t >  n\epsilon)\\
&\le& \sum_{t=1}^nP(Y_t> n\epsilon)
\;\le\; \sum \frac{EY^2_t}{n^2\epsilon^2}
\;=\; \frac{1}{n \epsilon^2}EY^2_t \longrightarrow 0
\end{eqnarray*}
as $EY^2_t$ =
   $E(X^4_{t-1}Z^4_t) < \infty$\\ \\
Conditions (4), (5), (6) and (7) follow from definition and condition on the
weights. This proves the claim.\\ \\
\noindent
Hence for $\tau^2 = \sigma^4(1-\theta^2)^{-1}$
\begin{equation} \label{eq1}
P\left[\;\frac{1}{\sqrt{n}}\sum_{t=2}^n W_t X_{t-1}Z_t \le x\;\vline\;X_1,\ldots, X_n\right] \stackrel{P}{\longrightarrow}\Phi(\;\;\frac{x}{\tau}\;\;) 
\hspace{.1in}\forall x \in \mathbb{R}
\end{equation}

\noindent
\textbf{Claim 2.} With $c = \sigma^2(1-\theta^2)^{-1}$,
\[P_B\left [\;\;\vline\;\frac{1}{n}\sum_{t=2}^{n} w_t X^2_{t-1} - c\;\vline\; > \;\epsilon\;\;\right ]
\stackrel{P}{\longrightarrow}0\hspace{.1in}\forall\; \epsilon > 0
\]
\textbf{Proof}
\begin{eqnarray*}
E_B(\frac{1}{n}\sum w_t X^2_{t-1}) &=& \frac{1}{n}\sum X^2_{t-1} \\
V_B(\sum w_t X^2_{t-1}) &=& \sum X^4_{t-1}\sigma^2_n + \sum _{s\;\neq\;t}X^2_{t-1}X^2_{s-1}\;Cov(w_t, w_s)\\
&=& \sigma^2_n \sum X^4_{t-1} + c_{1n}\sum_{s\;\neq\;t}X^2_{t-1}X^2_{s-1}
\end{eqnarray*}
Therefore 
\begin{eqnarray}
V_B(\frac{1}{n}\sum w_t X^2_{t-1}) &=& \frac{\sigma^2_n}{n^2}\sum X^4_{t-1}\hspace{.05in} + 
\hspace{.05in}\frac{c_{1n}}{n^2}\sum_{s\;\neq\;t}X^2_{t-1}X^2_{s-1} \label{veq}\\
\frac{1}{n}\sigma^2_n &\rightarrow& 0  \nonumber \\
\frac{1}{n}\sum X^4_{t-1} &\stackrel{a.s.}{\longrightarrow}& E(X^4_t) \nonumber
\end{eqnarray}
Hence the first term in (\ref{veq}) $\stackrel{a.s.}{\longrightarrow}0$\\
Also 
\begin{equation}
\frac{1}{n^2}\sum_{s\;\ne\;t}X^2_{t-1}X^2_{s-1}\; \le \; \left(\frac{\sum X^2_t}{n}\right)^2 \stackrel{a.s.}{\longrightarrow}(EX^2_t)^2 \nonumber
\end{equation}
Hence $\frac{1}{n^2}$ $\sum_{s\;\neq\;t}X^2_{t-1}X^2_{s-1}$ is bounded a.s., and
as $c_{1n}\rightarrow0$,

\smallskip
\noindent the second term in (\ref{veq}) also $\longrightarrow0$ a.s.

\smallskip
\noindent This shows that $V_B\left(\frac{1}{n}\sum w_t X^2_{t-1}\right)\longrightarrow0$ a.s.\\ \\
Hence $\frac{1}{n}\sum w_t X^2_{t-1}$ - $\frac{1}{n}\sum X^2_{t-1}\stackrel{P_B}{\longrightarrow}0$ a.s.\\ \\
Using Result (\ref{res2}), $\frac{1}{n}\sum X^2_{t-1}\stackrel{a.s.}{\longrightarrow}E(X^2_t)$ = $\sigma^2(1-\theta^2)^{-1}$ \\ \\
This implies, $\frac{1}{n}\sum w_t X^2_{t-1}$ $\stackrel{P_B}{\longrightarrow}\sigma^2(1-\theta^2)^{-1}$ a.s.

\smallskip 
\noindent This proves Claim 2.

\smallskip
\noindent In fact we have proved that, with $c = \sigma^2(1-\theta^2)^{-1}$
\begin{equation} \label{eq2}
P_B\left [\;\;\vline\; \frac{1}{n}\sum_{t=2}^{n} w_t X^2_{t-1} - c\; \vline\; > \; \epsilon\;\;\right ]
\stackrel{a.s.}{\longrightarrow}0\hspace{.1in}
\forall \; \epsilon > 0.
\end{equation}

\noindent
Now 
\begin{eqnarray}
& & \sqrt{n}\sigma^{-1}_n (\hat\theta^*_n-\hat{\theta}_n) \nonumber \\ 
&=& \sqrt{n}\sigma^{-1}_n\frac{\sum (w_t-1) X_{t-1} Z_t}{\sum w_t X^2_{t-1}}
- \sqrt{n}\sigma^{-1}_n \frac{\sum X_{t-1}Z_t \sum (w_t-1) X^2_{t-1}}{\sum  X^2_{t-1}\sum w_t X^2_{t-1}}\nonumber \\
&=& \frac{\sum W_t X_{t-1} Z_t/\sqrt{n}}{\sum w_t X^2_{t-1}/n}
 - \sqrt{n}(\hat{\theta}_n - \theta)\sigma^{-1}_n \frac{\sum (w_t-1) X^2_{t-1}/n}{\sum w_t X^2_{t-1}/n} \label{expand}\\
&=& T_1 - T_2\;\; (say) \nonumber
\end{eqnarray}
Then from (\ref{eq1}) and (\ref{eq2}), $ P_B(T_1 \le x) - P(T \le x) = o_P(1) $, where,
\begin{equation} \label{no}
T \sim \frac{1}{\sigma^2(1-\theta^2)^{-1}} N(0,\sigma^4(1-\theta^2)^{-1})\;=\;N(\;0,\;(1-\theta^2)\;)
\end{equation}
\textbf{Claim 3.} Define $A \equiv \sqrt{n}(\hat\theta_n - \theta)\sigma^{-1}_n\frac{1}{n}\sum(w_t -1)X^2_{t-1}$.\\

\smallskip
\hspace{.4in} Then $\forall \epsilon >0$,
$P_B(\;|A|>\epsilon \;)\stackrel{P}{\longrightarrow}0$.

\smallskip
\noindent \textbf{Proof} Note that,
\begin{eqnarray}
E_B(A) &=& 0 \\
V_B(A) &=& \frac{n}{\sigma^2_n}(\hat\theta_n - \theta)^2\;\; [ \frac{\sigma^2_n}{n^2}\sum X^4_{t-1} + \frac{c_{1n}}{n^2}\sum_s\sum_{\ne\;\; t} X^2_{s-1}X^2_{t-1} ] \\
&=& (\hat\theta_n - \theta)^2 \frac{\sum X^4_{t-1}}{n} + \frac{nc_{1n}}{\sigma^2_n}(\hat\theta_n - \theta)^2 \frac{\sum_s\sum_{\ne\;\; t} X^2_{s-1}X^2_{t-1}}{n^2} \\
&=& A_1 + A_2\;\;(say)
\end{eqnarray}
$\frac{\sum X^4_{t-1}}{n}$ converges a.s., and from (\ref{ars}), ($\hat\theta_n - \theta$)$\stackrel{P}{\longrightarrow}0$, as a result,
$A_1\stackrel{P}{\longrightarrow}0$.\\ \\
Moreover $\frac{\sum_{s\;\ne\;t} X^2_{s-1}X^2_{t-1}}{n^2}$ is bounded a.s., $nc_{1n}$ is bounded and $\sigma^2_n$ is
bounded away from 0. As a result $A_2\stackrel{P}{\longrightarrow}0$.

\smallskip
\noindent Combining, $V_B(A)\stackrel{P}{\longrightarrow}0$.

\smallskip
\noindent Hence 
\begin{eqnarray*}
P_B(\;|A|>\epsilon \;)\;\le\;\frac{V_B(A)}{\epsilon^2}\stackrel{P}{\longrightarrow}0
\end{eqnarray*}

\[\noindent Now \;\; T_2 = \frac{A}{\sum w_t X^2_{t-1}/n}.\]
From (\ref{eq2}), we have, $\sum w_t X^2_{t-1}/n$ is bounded away from zero in $P_B$ a.s., 
which means that, $\;\forall \epsilon > 0$,

\begin{equation} \label{no2}
P_B(\;|T_2| > \epsilon \;)\;=\;o_P(1)
\end{equation}
Hence from (\ref{expand}), (\ref{no}) and (\ref{no2}), we have,
\begin{equation} \label{bsres1}
\mathbf{
P_B[\;\sqrt{n}\sigma^{-1}_n(\hat\theta^{*}_n - \hat\theta_n) \le x\;] - P[\;Y \le x\;] = o_P(1)\hspace{.1in} \forall x \in \mathbb{R}
}
\end{equation}
where $Y \sim N(\;0,\;(1-\theta^2)\;)$ and this was what was to be proved.\\

\subsection{Least Absolute Deviations Estimator} Another estimator of $\theta_0$ can be the LAD estimatior, that is,
\[
\hat{\theta}_2 = \arg \min_{\theta}\frac{1}{n}\sum_2^n\;\vline X_t - \theta X_{t-1}\;\vline
\]
Now we reparametrize the model (\ref{ar}) in such a way that the median of $Z_t$, instead of the mean is equal to 0, while 
$VZ_t = \sigma^2$ remains unchanged.

\subsection {Distributional Consistency of the LAD estimator} Under the following assumptions we establish the assymptotic normality of  $\hat{\theta_2}$.
 \begin{description}
\item[A1.] CDF of $Z_t$, $F$ has a pdf $f$, which is continuous at zero. 
\item[A2.] $\vline\; F(x)\;-\; F(0) \; - \; xf(0) \; \vline\; \le\;  c |x|^{1+\alpha}$  in a neighborhood of zero, 
say $|x|\le M$, where $c \; , \; \alpha\; , M\; > 0$. 
\end{description}
\smallskip
To do so we use the the following result on random convex functions.

\begin{res}[See Niemire (1992)] \label{lem1}
Suppose that $h_n(a)$, $a \in R^d$ is a sequence of random convex functions which converge in probability to $h(a)$ for every fixed $a$. Then this convergence is uniform on any compact set containing $a$.
\end{res}
 
\begin{thm} \label{tl}
 Under the conditions (A1)-(A2), $\sqrt{n}(\hat{\theta}_2 \;-\; \theta_0) \stackrel{d}{\longrightarrow} N(0 , \frac{1}{4f^2(0)EX^2_t})$ as $n\rightarrow\infty$. 
\end{thm}
\noindent
\textbf{Proof}
Define
 \begin{eqnarray*}
 f(X_t , \theta) &=& ( |X_t - \theta X_{t-1}|\;-\;|X_t| ) \\
 g(X_t, \theta) &=& X_{t-1} [ 2I(Z_t(\theta) \le 0)\; -\; 1 ] \\
 where\; Z_t(\theta) &=& X_t - \theta X_{t-1}\;\; for t=2,\ldots,n. \\
 Y_t(a) &=&  f(X_t , \theta_0 + n^{-1/2}a)\;-\; f(X_t , \theta_0)\;-\;n^{-1/2}ag(X_t , \theta_0) \\
                       &=& |Z_t - n^{-1/2}aX_{t-1}|\;-\;|Z_t|\;-\;n^{-1/2}aX_{t-1} [ 2I(Z_t \le 0 ) - 1 ]\;\;for\; \alpha\; \in \mathbb{R}. 
\end{eqnarray*}
Also define
\begin{eqnarray*}
Q_n(\theta) &=& \sum f( X_t,\theta ) \\
U_n &=& \sum g(X_t , \theta_0) \\
V_n &=& \sum Y_t(a) = Q_n( \theta_0 + n^{-1/2}a )\;-\;Q_n(\theta_0)\;-\;n^{-1/2}a U_n
\end{eqnarray*}
  
\noindent \textbf{Step1} $\sum_{t=2}^n Y_t(a)\stackrel{P}{\longrightarrow}a^2f(0)EX^2_1$ \\
\smallskip
\textbf{Step1.1} $\sum ( Y_t\;-\;E(Y_t|\mathcal{A}_{t-1}) ) \stackrel{P}{\longrightarrow}0$
 \begin{eqnarray*}
 E(\;Y_t\;-\;E(Y_t|\mathcal{A}_{t-1}) \;) &=& 0 \\
 V(\;\sum Y_t\;-\;E(Y_t|\mathcal{A}_{t-1}) \;) &=& \sum V(\;Y_t\;-\;E(Y_t|\mathcal{A}_{t-1}) \;)\;\le\;\sum V(Y_t)\;\le\;\sum EY^2_t 
 \end{eqnarray*}
 By convexity of $f$, 
 \[ 0 \le Y_t(a) \le n^{-1/2}a [ g(X_t , \theta_0 + n^{-1/2}a ) - g(X_t , \theta_0 )] \]
 Therefore 
 \begin{eqnarray*}
  E(Y^2_t) &\le& \frac{ a^2 }{n} E[  g(X_t , \theta_0 + n^{-1/2}a ) - g(X_t , \theta_0 ) ]^2 \\
                   &=& 4\frac{a^2 }{n} E  X^2_{t-1} [ I(Z_t - n^{-1/2}a X_{t-1}\le 0)\;-\;I(Z_t \le 0) ]^2
\end{eqnarray*}
Now
\[ \sum EY^2_t = nEY^2_2 \le  4 a E  X^2_1 [ I(Z_2 - n^{-1/2}a X_1\le 0)\;-\;I(Z_2 \le 0) ]^2 \]
which tends to zero using DCT.
Therefore \[ V(\;\sum ( Y_t\;-\;E(Y_t|\mathcal{A}_{t-1}) \;) \rightarrow0 \]
This establishes Step 1.1. 

\smallskip
\noindent
\textbf{Step1.2} $\sum E(Y_t\vline \mathcal{A}_{t-1}) - a^2f(0)EX^2_1 \stackrel{P}{\longrightarrow}0$ 
\begin{eqnarray*}
E(Y_t|\mathcal{A}_{t-1}) &=& E(|Z_t - n^{-1/2}a X_{t-1}|\vline \mathcal{A}_{t-1}) - E|Z_t| \\
&=& \int ( |z - n^{-1/2}a X_{t-1}| - |z|) dF(z)
\end{eqnarray*}
Using the representation,
\[
| x-\theta | -|x| = \theta [ 2I(x \le 0) - 1 ] + 2\int_0^{\theta}[ I(x \le s) - I(x \le 0) ]ds
\]
we have
\[
|z - n^{-1/2}a X_{t-1}| - |z| = n^{-1/2}a X_{t-1} [ 2I(z \le 0) - 1 ] + 2\int_0^{n^{-1/2}a X_{t-1}}[ I(z \le s) - I(z \le 0)]ds
\] 
Therefore
\begin{eqnarray}
E(Y_t|\mathcal{A}_{t-1}) &=& n^{-1/2}a X_{t-1}\int [ 2I(z \le 0) - 1 ]dFz \;+\;  2\int
\int_0^{n^{-1/2}a X_{t-1}}[ I(z \le s) - I(z \le 0)]dsdFz \\ 
&=& 2\int_0^{n^{-1/2}a X_{t-1}}[F(s) - F(0)]ds \\
&=& 2n^{-1/2}X_{t-1}\int_0^a[F(n^{-1/2}X_{t-1}x) - F(0)]dx \label{bsn}
\end{eqnarray}
Under assumption A2, 
\begin{eqnarray*}
F(n^{-1/2}X_{t-1}x) - F(0) &=& n^{-1/2}X_{t-1}xf(0) + R_{nt}(x) \\
where\;|R_{nt}(x)| &\le& cn^{-(1+\alpha)/2}|X_{t-1}|^{1+\alpha}|x|^{1+\alpha} \\
whenever \; n^{-1/2}|X_{t-1}||x| &\le& M
\end{eqnarray*}
Hence
\begin{eqnarray*}  
E(Y_t|\mathcal{A}_{t-1}) &=& 2n^{-1/2}X_{t-1}\int_0^a[n^{-1/2}X_{t-1}xf(0) + R_{nt}(x)]dx \\
&=& \frac{1}{n}X^2_{t-1}a^2f(0) + 2n^{-1/2}X_{t-1}\int_0^aR_{nt}(x)dx \\
\sum E(Y_t|\mathcal{A}_{t-1}) &=& a^2f(0)\frac{1}{n}\sum X^2_{t-1} + \frac{2}{n}\sum X_{t-1}\int_0^a\sqrt{n}R_{nt}(x)dx \\
&=& I_1 + I_2 \;(say)
\end{eqnarray*}
Then $I_1\stackrel{P}{\longrightarrow}a^2f(0)EX^2_1$.\\
Remains to show $I_2\stackrel{P}{\longrightarrow}0$.
To show this, let us assume:
\begin{description}
\item[1.] $\max_{1\le t\le n}\;n^{-1/2}|X_{t-1}|\stackrel{P}{\longrightarrow}0$ 
\item[2.] $\frac{1}{n^{1+\alpha/2}}\sum |X_{t-1}|^{2+\alpha}\stackrel{P}{\longrightarrow}0$ 
\end{description}
Hence given $\epsilon>0$, 
\begin{eqnarray*}
P(\max\; n^{-1/2}|X_{t-1}|\le M/|a|)&\rightarrow & 1 \\
and \; P( \frac{c}{n^{1+\alpha/2}}\sum |X_{t-1}|^{2+\alpha} < \epsilon ) &\rightarrow& 1 
\end{eqnarray*}
Let $A_n$ be the set where $max\; n^{-1/2}|X_{t-1}|\le M/|a|$ and 
$\frac{c}{n^{1+\alpha/2}}\sum |X_{t-1}|^{2+\alpha} < \epsilon$.\\
Then $\exists N$ such that $P(A_n)> 1-\epsilon\; \forall \; n \ge N$.
Then on $A_n$, $|R_{nt}| \le cn^{-\alpha/2}|X_{t-1}|^{1+\alpha}$, and hence
\begin{eqnarray*}
|I_2| &\le& \frac{2}{n} \sum |X_{t-1}|\int_0^a cn^{-\alpha/2}|X_{t-1}|^{1+\alpha}\\
&\le& \frac{c}{n^{1+\alpha/2}}\sum |X_{t-1}|^{2+\alpha} \\
&<& \epsilon
\end{eqnarray*}
ie $P( |I_2| < \epsilon ) \rightarrow1 \; \forall \; \epsilon > 0$. In otherwords $I_2 \stackrel{P}{\longrightarrow}0$.
This completes Step1.2 and hence Step1.
\smallskip 
In other words, 
\begin{equation}
Q_n( \theta_0 + n^{-1/2}a )\;-\;Q_n(\theta_0)\;-\;n^{-1/2}a U_n - a^2f(0)EX^2_1 \stackrel{P}{\longrightarrow}0  \label{ow}
\end{equation}
Due to convexity of $Q_n$, the convergence in (\ref{ow}) is uniform on any compact set by Result\ref{lem1}. Thus $\forall$ $\epsilon > 0$, and $M > 0$, for $n$ sufficiently large, we have
\[
P \left[ \sup_{|a| \le M} \left| Q_n( \theta_0 + n^{-1/2}a )\;-\;Q_n(\theta_0)\;-\;n^{-1/2}a U_n - a^2f(0)EX^2_1 \right| < \epsilon \right] \ge 1 - \epsilon /2
\]
Call
\begin{eqnarray*}
A_n(a) &=& Q_n( \theta_0 + n^{-1/2}a )\;-\;Q_n(\theta_0)\;\; ,\\
B_n(a) &=& n^{-1/2}a U_n + a^2f(0)EX^2_1
\end{eqnarray*}
and their minimizers $a_n$ and $b_n$ respectively. Then
\begin{eqnarray*}
a_n &=& \sqrt{n} (\hat\theta_2 - \theta_0) \;and \\
b_n &=& -(2f(0)EX^2_1)^{-1} n^{-1/2} U_n
\end{eqnarray*}
The minimum value of $B_n$,
\[ B_n(b_n) = -n^{-1} (4f(0)EX^2_1)^{-1}U^2_n \]
Note that $b_n$ is bounded in probability. Hence there exists $M > 0$ such that
\[
P \left[ |-(2f(0)EX^2_1)^{-1} n^{-1/2} U_n| < M-1 \right] \ge 1- \epsilon /2 \]
Let $A$ be the set where,
  
\begin{eqnarray*}
&& \sup_{|a|\le M} |A_n(a) - B_n(a)| < \epsilon \\
and\; && |-(2f(0)EX^2_1)^{-1} n^{-1/2} U_n| < M-1
\end{eqnarray*}
Then $P(A) > 1-\epsilon$.
On $A$, 
\begin{equation}
A_n(b_n) < B_n(b_n) + \epsilon \label{ab1}
\end{equation}
Consider the value of $A_n$ on the sphere $ S_n = \{a : |a - b_n| = k\epsilon^{1/2}\}$ where $k$ will be chosen later. By chosing  $\epsilon$ sufficiently small, we have $|a| \le M$ $\forall a \in S_n$. Hence 
\begin{equation}
A_n(a) > B_n(a) - \epsilon\; \forall a \in S_n. \label{ab2}
\end{equation}
Once we chose $k = 2(2f(0)EX^2_1)^{-1/2}$, 
\begin{equation}
B_n(a) > B_n(b_n) + 2\epsilon \; \forall a \; \in S_n
\end{equation}
Comparing the bounds (\ref{ab1}) and (\ref{ab2}), we have $A_n(a) > A_n(b_n)$ whenever $a \in S_n$. If $|a_n - b_n| > k\epsilon^{1/2}$, by convexity of $A_n$, thereexists $a^{*}_n$ on $S_n$ such that $A_n(a^{*}_n) \le A_n(b_n)$ which cannot be the case. Therefore $|a_n - b_n| < k\epsilon^{1/2}$ on $A$.
Since this holds with probability atleast $1 - \epsilon$ and $\epsilon$ is arbitrary, this proves that $|a_n - b_n| \stackrel{P}{\longrightarrow} 0$. In otherwords,
 \begin{eqnarray}
 \sqrt{n}(\hat{\theta_2} - \theta_0) = -n^{-1/2}(2f(0)EX^2_1)^{-1}U_n + o_P(1) \label{el1}
 \end{eqnarray}
 \smallskip
 \textbf{Step 2} $n^{-1/2}U_n \stackrel{d}{\longrightarrow}N(0 , EX^2_1)$ 
 \begin{eqnarray*}
 U_n &=& \sum_{t=2}^nX_{t-1}[2I(Z_t \le 0) - 1] \\
 &=& \sum_{t=2}^nY_t \;(say)
 \end{eqnarray*}
 Then note that $U_n$ is a 0-mean martingale with finite variance increments. Hence to prove Step2, we use the Martingale CLT.
 Write 
 \begin{eqnarray*}
 S^2_n &=& \sum_{t=2}^nE(Y^2_t|\mathcal{A}_{t-1}) = \sum_{t=2}^nX^2_{t-1} \\
 and \; s^2_n &=& ES^2_n = (n-1)EX^2_1
 \end{eqnarray*}
 Then we need to to verify:
 \begin{description}
 \item[1.] $\frac{S^2_n}{s^2_n} \stackrel{P}{\longrightarrow}1$ \\
 \noindent 
 This follows from Result \ref{res2}.
 \item[2.] $s^{-2}_n\sum_{t=2}^{n}E(Y^2_tI(|Y_t|\ge \epsilon s_n)) \longrightarrow 0\;\; as \; n
 \rightarrow \infty \;\;\forall \epsilon > 0$. \\
 \noindent 
 To see this, note that 
 \begin{eqnarray*}
 L.H.S. \; &=& \frac{1}{EX^2_1}E(X^2_1 I\left(\frac{|X_1|}{\sqrt{EX^2_1}} \ge \epsilon \sqrt{n-1} \right)) \\ 
 &\longrightarrow&0 \; as \; EX^2_1 < \infty
 \end{eqnarray*}
 \end{description} 
 Hence using Result \ref{mclt}, we have $ \frac{U_n}{s_n} \stackrel{d}{\longrightarrow}N(0,1) $, which proves Step2.\\
 Combining Step2 and equation(\ref{el1}), we get,
 \[
 \sqrt{n}(\hat{\theta_2} - \theta_0) \stackrel{d}{\longrightarrow} N\left(0 , \frac{1}{4f^2(0)EX^2_1}\right)
 \]
 and this was what was to be proved. \\
 \medskip
 Finally it remains to verify: \\
 $1. \; \max_{2\le t \le n} n^{-1/2} |X_{t-1}|\stackrel{P}{\longrightarrow}0$ \\
 Proof: Given $\epsilon$ positive,
 \begin{eqnarray*}
 P( \max_t n^{-1/2} |X_{t-1} > \epsilon ) &\le& \sum_{t=1}^{n-1} 	P(|X_t| > \epsilon\sqrt{n}) \\
 &=& (n-1)P(|X_1| > \epsilon\sqrt{n}) \\
 &=& (n-1) \int I(|X_1| > \epsilon \sqrt{n})dP \\
 &\le& (n-1) \int \frac{|X_1|^2}{\epsilon^2 n} I(|X_1| > \epsilon \sqrt{n} dP \\
 &=& \frac{1}{\epsilon^2} \int |X_1^2 I( |X_1| > \epsilon \sqrt{n} ) dP \\
 &\rightarrow& 0 \; as \; E|X_1|^2 < \infty
 \end{eqnarray*}
 \medskip
 \noindent
 $2. \; \frac{1}{n^{1+\alpha/2}} \sum_{t=2}^n |X_{t-1}|^{2+\alpha} \stackrel{P}{\longrightarrow}0$\\
 Proof:                                      
\begin{eqnarray*}
\frac{1}{n^{1+\alpha/2}} \sum_{t=2}^n |X_{t-1}|^{2+\alpha} &\le& \frac{\max_{1\le t \le n-1} |X_t|^{\alpha}}{n^{1+\alpha/2}} \frac{1}{n} \sum X_{t-1}^2 \\
&\le& \left( \frac{\max |X_t|}{\sqrt{n}} \right)^{\alpha} \frac{1}{n} \sum X_{t-1}^2  \\
&\stackrel{P}{\longrightarrow}0
\end{eqnarray*}
This follows from (1) and the fact that $\frac{1}{n} \sum X_{t-1}^2$ is bounded in probability, since 
$EX^2_1 < \infty$.
This completes the proof.

\subsection{WBS for LAD estimators}
Now we define the weighted bootstrap estimators, $\hat\theta^{*}_2$ of $\hat\theta_2$ as the minimizers of
\begin{equation}\label{wbl}
Q_{nB}(\theta) = \sum_{t=2}^n w_{nt}|X_t - \theta X_{t-1}|
\end{equation}  
In the next section, we deduce the consistency of this bootstrap procedure.

\subsection{ Consistency of the Weighted Bootstrap technique}
Now we prove that the Weighted Bootstrap estimator of $\hat{\theta_2}$ is assymptotically normal with the same assymptotic
distribution. In particular WB provides a consistent resampling scheme to estimate the LAD estimator.

\noindent
\begin{thm} \label{tbl}
Let $\hat\theta^{*}_2$ be the weighted bootstrap estimator of $\hat\theta_2$ as defined in (\ref{wbl}). Suppose the bootstrap weights 
satisfy conditions (A1)-(A4). Also assume that $n^{-1/2}\sigma_n \max_t |X_t| \stackrel{P}{\longrightarrow}0$. Then
\begin{equation}
\sup_{x \in \mathbb{R}}\;\vline P\left [\;\sqrt{n}\sigma^{-1}_n(\hat\theta^{*}_2 - \hat\theta_2) \le x\;
\vline X_1,\ldots,X_n\right ] - P\left [Y \le x\right ]\; \vline =
 o_P(1)
\end{equation}
where $Y \sim N \left(0, \frac{1}{4f^2(0)EX^2_t} \right)$.
\end{thm}

\noindent
\textbf{Proof} Define
\begin{eqnarray*}
U_{nt}(a) &=& f(X_t , \theta_0 +n^{-1/2}\sigma_n a) - f(X_t , \theta_0) - n^{-1/2} \sigma_n a g(X_t , \theta_0) \\
U_{nBt}(a) &=& w_{nt} U_{nt}(a) \\
S_{nB} &=& \sum W_{nt} g(X_t , \theta_0) \\
S_{nw} &=& \sum w_{nt} g(X_t , \theta_0) \\
S_n &=& \sum g(X_t , \theta_0) \\
H &=& 2 f(0)EX^2_1
\end{eqnarray*}

Then
\begin{eqnarray*}
E_B U_{nBt} &=& U_{nt} \; and \\
\sum U_{nBt}(a) &=& Q_{nB}(\theta_0 + n^{-1/2}\sigma_n a) - Q_{nB}(\theta_0) - n^{-1/2} \sigma_n aS_{nw}
\end{eqnarray*}  

\noindent
\textbf{Step1.} We show $\sqrt{n}\sigma^{-1}_n(\hat\theta^{*}_2 - \hat\theta_2) = -n^{-1/2} H^{-1} S_{nB} + r_{nB}$ s.t. given   
$\epsilon > 0$, $P_B [ |r_{nB}|> \epsilon ] = o_P(1)$.\\
\smallskip
\noindent
To show this, choose $k = 3H^{-1/2}$ and $\epsilon$ small enough such that $k^2\epsilon < 1$ and $M$ a sufficiently large constant. Let 
$\mathcal{A}$ be the set where
\begin{eqnarray*}
&& \sup_{|a| \le M} \sigma^2_n \left| Q_{nB}(\theta_0 + n^{-1/2}\sigma_n a) - Q_{nB}(\theta_0) - n^{-1/2}\sigma_n a 
S_{nw} - \frac{\sigma^2_n}{2}a^2H \right| < \epsilon \\
&and& \left| n^{-1/2} \sigma^{-1}_n H^{-1} S_{nw} \right| < M-1
\end{eqnarray*} 

Then due to convexity of $Q_{nB}$, arguing as in the proof of Theorem2,on $\mathcal{A}$ we have,
\begin{eqnarray*}
&& \sqrt{n}\sigma^{-1}_n (\hat\theta^{*}_2 - \theta_0) = -n^{-1/2} \sigma^{-1}_nH^{-1}S_{nW} + r_{nB} \\
&s.t.& |r_{nB}| < k \epsilon^{1/2}
\end{eqnarray*}

If we show $1 - P_B[ \mathcal{A} ] = o_P(1)$, then 
\[
P_B[ |r_{nB}| > \delta] = o_P(1) \; \; \forall \; \delta > 0 
\]
Also from equation(\ref{el1}); $\sqrt{n}( \hat\theta_2 - \theta_0 ) = -n^{-1/2}H^{-1} S_n + o_P(1)$.\\
Therefore $\sqrt{n}\sigma^{-1}_n(\hat\theta^{*}_2 - \hat\theta_2) = -n^{-1/2} H^{-1} S_{nB} + r_{nB2}$ s.t. given
$\epsilon > 0$, $P_B [ |r_{nB2}|> \epsilon ] = o_P(1)$.\\
This will complete Step1.\\ \\
Hence it remains to show, $1 - P_B[\mathcal{A}] = o_P(1)$\\
To show this we show,
\begin{eqnarray}
&\forall& M > 0 , \; P_B\left[sup_{|a|\le M} \sigma^{-2}_n \left|\sum U_{nBt}(a) - \frac{\sigma^2_n}{2}a^2H\right| > \epsilon \right] = o_P(1) \label{ts1} \\
&and& \; there\;exists\; M > 0 \; s.t.\; P_B \left[ |\sigma^{-1}_n n^{-1/2} H^{-1}S_{nw}| \ge M \right] = o_P(1) \label{ts2}
\end{eqnarray}
To show (\ref{ts1}), note that,
\begin{eqnarray*}
&& P_B\left[sup_{|a|\le M} \sigma^{-2}_n \left|\sum U_{nBt}(a) - \frac{\sigma^2_n}{2}a^2H\right| > \epsilon \right] \\
&\le& \sum_j P_B[\sigma^{-1}_n |\sum_t W_t U_t(b_j)| > \epsilon/2] + \sum_j I( \sigma^{-2}_n |\sum X_t(b_j) - \sigma^2_n b^2_j H/2 | > \epsilon/2 ) \\
&\le& \sigma^{-2}_n \sum_j k \sum_t U^2_t(b_j) + \sum_j I( \sigma^{-2}_n |\sum_t U_t(b_j) - \sigma^2_n b^2_j H/2 | > \epsilon/2) 
\end{eqnarray*}
As a result, we need to show for fixed $b$,
\begin{eqnarray}
\sigma^{-2}_n \sum_t U^2_{nt}(b) &=& o_P(1) \label{cl1} \\
and \; \sigma^{-2}_n [ \sum_t U_{nt}(b) - \sigma^2_n b^2 H/2 ] &=& o_P(1) \label{cl2}
\end{eqnarray}
 
To see (\ref{cl1}), 
\begin{eqnarray*}
\sigma^{-2}_n \sum_t EU^2_{nt}(b) &=& n \sigma^{-2}_n EU^2_1(b) \\
&\le& n \sigma^{-2}_n E [ f(X_1 , \theta_0 + n^{-1/2}\sigma_n b) - f(X_1 , \theta_0) - n^{-1/2} \sigma_n b 
g(X_1, \theta_0) ]^2 \\
&\le& E b^2 [ g(X_1 , \theta_0 + n^{-1/2}\sigma_n b) - g(X_t , \theta_0) ]^2 \\
&\longrightarrow& 0 
\end{eqnarray*}
This proves (\ref{cl1}). \\ \\

\noindent To prove (\ref{cl2}) note that,

\begin{eqnarray*}
&& \sigma^{-2}_n \left[ \sum U_t(b) - \sigma^2_n b^2 H/2 \right] = \sigma^{-2}_n \left[ \sum [U_t(b) - E(U_t(b)|\mathcal{A}_{t-1})] + \sum E( U_t(b)|\mathcal{A}_{t-1}) - \sigma^2_n b^2 H/2 \right] \\
&& E \sigma^{-2}_n \sum (U_t(b) - E(U_t(b)|\mathcal{A}_{t-1})) = 0 \\
&& V[ \sigma^{-2}_n \sum (U_t(b) - E(U_t(b)|\mathcal{A}_{t-1})) ] = \sigma^{-4}_n \sum V(U_t - E(U_t|\mathcal{A}_{t-1})) \\
&\le& \sigma^{-4}_n \sum V(U_t(b)) \\
&\le& k^{-1}_1 \sigma^{-2}_n \sum E(U^2_t(b)) \;\; \; (\sigma^2_n > k_1) \\
&=& nk^{-1}_1 \sigma^{-2}_n E(U_1(b))^2 \\
&\le& k^{-1}_1 \sigma^{-2}_n \sigma^2_n b^2 E[ g(X_1, \theta_0 + n^{-1/2}\sigma_n b) - g(X_1, \theta_0)]^2 \\
&=& \frac{1}{k_1} E[ g(X_1, \theta_0 + n^{-1/2}\sigma_n b) - g(X_1, \theta_0) ]^2 \\
&\longrightarrow& 0
\end{eqnarray*}
Hence
\begin{equation}
\sigma^{-2}_n [ \sum( U_1(b) - E( U_t(b)|\mathcal{A}_{t-1}) ) ] \stackrel{P}{\longrightarrow}0 \label{cl21}
\end{equation}

\begin{eqnarray}
&& \sigma^{-2}_n\sum E(U_t|\mathcal{A}_{t-1}) = \sigma^{-2}_n \sum 2n^{-1/2}X_{t-1} \int_0^{\sigma_n b} [ F(n^{-1/2}\sigma_n X_{t-1}x) - F(0) ]dx \label{c221} \\
&=& 2n^{-1/2}\sigma^{-1}_n \sum X_{t-1} \int_0^b [ F(n^{-1/2}\sigma_n X_{t-1}x) - F(0) ]dx \\ 
&=& 2n^{-1/2}\sigma^{-1}_n \sum X_{t-1} \left[ n^{-1/2}\sigma_n \frac{b^2}{2}X_{t-1}f(0) + R_{nt} \right] \\
&where& \; |R_{nt}| \le c|n^{-1/2}\sigma_n X_{t-1}|^{1+\alpha} = c (n^{-1/2}\sigma_n)^{1+\alpha} |X_{t-1}|^{1+\alpha} 
\label{c222} \\
&=& 2 \frac{b^2}{2} f(0)\frac{1}{n} \sum X^2_{t-1} + 2 n^{-1/2} \sigma^{-1}_n \sum X_{t-1} R_{nt} \\
&=& I_1 + I_2 \;\;(say)
\end{eqnarray}
Here (\ref{c221}) follows from (\ref{bsn}), and (\ref{c222}) from assumption A2 on $F$ and the assumption 
$n^{-1/2}\sigma_n \max_t |X_t| \stackrel{P}{\longrightarrow}0$.

\begin{eqnarray}
Now\;\; I_1 &\stackrel{P}{\longrightarrow}& b^2f(0)EX^2_1 = b^2H/2 \label{i11} \\
and \;\; |I_2| &\le& c (n^{-1/2}\sigma^{-1}_n)(n^{-1/2}\sigma_n)^{1+\alpha} \sum |X_{t-1}|^{2+\alpha} \\
&=& \frac{c \sigma^{\alpha}_n}{n^{1 + \alpha/2}} \sum |X_{t-1}|^{2+\alpha} \label{i2} \\
&\stackrel{P}{\longrightarrow}&0 \label{i22}
\end{eqnarray}

\noindent In this case (\ref{i22}) follows from (\ref{i2}) if we show $\frac{\sigma^{\alpha}_n}{n^{1 + \alpha/2}} \sum |X_{t-1}|^{2+\alpha} \stackrel{P}{\longrightarrow}0$ $\forall \alpha>0$.\\
To see this, note that
\begin{eqnarray*}
\frac{\sigma^{\alpha}_n}{n^{1 + \alpha/2}} \sum |X_{t-1}|^{2+\alpha} &\le& \frac{\sigma^{\alpha}_n \max_t |X_t|^{\alpha}}{n^{\alpha/2}} \frac{1}{n}\sum X^2_t \\ 
&\le& \left(\frac{\sigma_n}{\sqrt{n}} \max_{1\le t\le n}|X_t|\right)^{\alpha} \frac{1}{n} \sum X^2_t \\
&\stackrel{P}{\longrightarrow}0
\end{eqnarray*}
Combining (\ref{i11}) and (\ref{i22}), we have $\sigma^{-2}_n\sum E(U_t|\mathcal{A}_{t-1}) \stackrel{P}{\longrightarrow} b^2 H/2$. In o.w.,
\begin{equation}
\sigma^{-2}_n \left[ \sum E(U_t |\mathcal{A}_{t-1}) - \sigma^2_n b^2 H/2 \right] \stackrel{P}{\longrightarrow}0 \label{cl22}
\end{equation}
Adding (\ref{cl21}) and (\ref{cl22}), we prove (\ref{cl2}). And from (\ref{cl1}) and (\ref{cl2}) we deduce (\ref{ts1}).

\begin{eqnarray*}
&Now& \; \; P_B\left[ |\sigma^{-1}_n n^{-1/2} H^{-1} S_{nw}| \ge M \right] \\
&\le& \frac{\sigma^{-2}_n n^{-1} H^{-2}}{M^2} E_B \left[ \sum w_t g(X_t, \theta_0) \right]^2 \\
&\le& \frac{K_1}{M^2 n} \sum g(X_t , \theta_0)^2 + \frac{K_2}{M^2} \left[ \frac{\sum g(X_t , \theta_0)}{\sqrt{n}} \right]^2 
\\
&\stackrel{P}{\longrightarrow}0
\end{eqnarray*}
if $M$ is choosen sufficiently large. This proves (\ref{ts2}). (\ref{ts1}) and (\ref{ts2}) together show $1 - P_B[\mathcal{A}] = o_P(1)$. This completes step 1. \\ \\

\noindent \textbf{Step2.} $ P_B( n^{-1/2}S_{nB} \le x ) - P(Y \le x ) = o_P(1)$, where $Y \sim N(0, EX^2_1)$ \\
\smallskip
\noindent To show this we use Result1.
\begin{eqnarray*}
S_{nB} &=& \sum W_{nt} g(X_t, \theta_0)  \\
&=& \sum W_{nt} X_{t-1} [ 2I(Z_t \le 0) - 1] 
\end{eqnarray*}
Hence we need to show:
\begin{enumerate}
\item $\frac{1}{n}\sum X_{t-1}[ 2I(Z_t \le 0) - 1]\stackrel{P}{\longrightarrow}0$
\item $\frac{1}{n}\sum X^2_{t-1} \stackrel{P}{\longrightarrow}EX^2_1$
\item $\frac{1}{n} \max_t X^2_{t-1}\stackrel{P}{\longrightarrow}0$
\end{enumerate}
All these follow from Step2 in the proof of Theorem2. \\ \\
This completes Step2, and combining with Step1, we get
\begin{eqnarray*}
&& P_B \left(\sqrt{n}\sigma^{-1}_n (\hat\theta^{*}_2 - \hat\theta_2) \le x \right) - P(Y \le x) = o_P(1)\; \forall \; x \; \in \; \mathbb{R}\\
where && \; Y \sim N\left(0,\frac{1}{4f^2(0)EX^2_1}\right)
\end{eqnarray*}
Using continuity of the normal distribution, we complete the proof.

\subsection{Special choices for w.} With $(w_1,\ldots,w_n) \sim Mult(n,\frac{1}{n},\ldots,\frac{1}{n})$ we get the Paired Bootstrap estimator. This is same as resampling w.r. from ($X_{t-1}$, $X_t$), 
$t = 1,2,\ldots,n$. Other choices of $\{w_i\}$'s yield the m-out-of-n Bootstrap and their variations.
In particular lets check the conditions on the weights in two particular cases.

\medskip
\noindent
\textbf{Case 1.\hspace{0.1in}$\mathbf{(w_1,\ldots,w_n) \sim Mult(n,\frac{1}{n},\ldots,\frac{1}{n})}$}

\smallskip
\noindent Clearly the weights are exchangeable. Let us verify assumptions (A1)-(A4) on the weights in this case.
\begin{description}
\item[A1.] $E_B(w_1)=1$\\
Obvious in this case.
\item[A2.] $0 < k < \sigma^2_n = o(n)$\\
$\sigma^2_n = 1-\frac{1}{n}$ which clearly satisfies the above condition.
\item[A3.] $c_{1n} = O(n^{-1})$\\
$c_{1n} = -\frac{1}{n}$ which is as above.
\item[A4.] $\{W_i\}$ satisfy conditions of P-W theorem.\\
To show this, we have to verify conditions (6) and (7) of Result (\ref{pw}) with $U_{nj}\;=\; W_j$.\\ \\ 
\textbf{Condition(6)} $\frac{1}{n}\sum{W^2_t}{\stackrel{P}{\longrightarrow}}1$\\
$W_t = \sqrt{\frac{n}{n-1}}(w_t-1)$\\
Therefore,\[\frac{1}{n}\sum W^2_t = \frac{n}{n-1}\frac{1}{n}\sum(w_t-1)^2\]
\[V_B(\sum(w_t-1)^2) = nV_B(w_1-1)^2 + n(n-1)Cov_B((w_1-1)^2,(w_2-1)^2)\]
Write $w_1 = \sum_{i=1}^n u_i$ and $w_2 = \sum_{i=1}^n v_i$, where $\{u_i,v_i\}_{i=1}^n$ are iid with the joint 
distribution of ($u_i,v_i$) given by
\[
  (u_i,v_i) = 
  \left\{
  \begin{array}{ll}
   (1,0) & \mbox{w.p. $1/n$}\\
   (0,1) & \mbox{w.p. $1/n$}\\
   (1,0) & \mbox{w.p. $1 - 2/n$}
   \end{array}
   \right.
 \]
\begin{eqnarray*}
V_B(w_1-1)^2 &=& E_B(w_1-1)^4 - V^2_B(w_1) \\
&=& E_B(w_1-1)^4 - \frac{1}{n^2} (1-\frac{1}{n})^2 
\end{eqnarray*}
$(w_1-1) = \sum_{i=1}^n (u_i - p)$ where $p\;=\;\frac{1}{n}$ , $q \;=\; 1 - p$. Hence
\begin{eqnarray*}
E_B(w_1-1)^4 &=& E(\;\; \sum(u_i-p)^4 + 3\sum_{i\;\ne \; j}(u_i-p)^2(u_j-p)^2 \;\;) \\
&=& nE(u_1-p)^4 + 3n(n-1)p^2q^2 \\
&=& n(pq^4 + p^4q) + 3n(n-1)p^2q^2
\end{eqnarray*}

\noindent Simplifying
\begin{equation} \label{mom4}
= (1-\frac{1}{n})(4 - \frac{9}{n} + \frac{6}{n^2} + \frac{2}{n^3}) 
\end{equation}
Therefore
\begin{equation}
V_B(w_1-1)^2 = (1-\frac{1}{n})(3 - \frac{9}{n} + \frac{7}{n^2} + \frac{2}{n^3})\rightarrow3
\end{equation}

\begin{eqnarray*}
(w_1-1)^2(w_2-1)^2 &=& [\sum(u_i-p)]^2[\sum(v_i-p)]^2 \\
&=& [\sum(u_i-p)^2 + \sum_{i\ne j}\sum(u_i-p)(u_j-p)] \\
{} &\times& [\sum(v_i-p)^2 + \sum_{i\ne j}\sum(v_i-p)(v_j-p)]
\end{eqnarray*}

\begin{eqnarray*}
& & E_B(w_1-1)^2(w_2-1)^2 \\ 
&=& E[\;\sum_{i}(u_i-p)^2(v_i-p)^2 \;+\; \sum_{\;i\;\ne}\sum_{\;j}(u_i-p)^2(v_j-p)^2 \\
{}& & +\; \sum_{\;i\;\ne}\sum_{\;j}(u_i-p)(u_j-p)(v_i-p)(v_j-p)\;] \\
&=& nE(u_1-p)^2(v_1-p)^2 \;+\; n(n-1)V(u_1)V(v_1) \\
{}& & +\; n(n-1)Cov(u_1,v_1)Cov(u_2,v_2) \\
&=& n(2pp^2q^2 \; + \; p^4(1-2p))\; + \;n(n-1)p^2q^2 \;-\; n(n-1)p^4 \\
&=& \frac{2}{n^2}(1-\frac{1}{n})^2 \;+\;\frac{1}{n^3}(1-\frac{2}{n})\\
{}& & +\; n(n-1)\frac{1}{n^2}(1-\frac{1}{n})^2\;- \;n(n-1)\frac{1}{n^4} \\
&=& 1 \;-\; \frac{3}{n} \;+\; \frac{4}{n^2} \;-\; \frac{3}{n^3}
\end{eqnarray*} 

\begin{eqnarray*}
& & Cov_B((w_1-1)^2,(w_2-1)^2) \\
&=& E_B(w_1-1)^2(w_2-1)^2 \;-\; (1-\frac{1}{n})^2 \\
&=& 1 \;-\; \frac{3}{n} \;+\; \frac{4}{n^2} \;-\; \frac{3}{n^3} \\
{}& & -\; (1- \frac{2}{n} \;+\; \frac{1}{n^2}) \\
&=& -\;\frac{1}{n} \;+\; \frac{3}{n^2} \;-\; \frac{3}{n^3} \\
&\longrightarrow& 0
\end{eqnarray*} 
Therefore
\begin{eqnarray}
V_B(\frac{1}{n}\sum(w_t-1)^2) &=& \frac{1}{n}V_B(w_1-1)^2 +
(1-\frac{1}{n})Cov_B((w_1-1)^2,(w_2-1)^2) \nonumber \\
& & \longrightarrow0 \nonumber \\
V_B(\frac{1}{n}\sum W^2_t) &=& 
(\frac{n}{n-1})^2 \; V_B(\frac{1}{n}\sum(w_t-1)^2)\longrightarrow0
\label{var}\\
E_B(\frac{1}{n}\sum W^2_t) &=& \frac{n}{n-1}E_B(w_1-1)^2 = 1 \label{exp}
\end{eqnarray}

\noindent
Hence from (\ref{var}) and (\ref{exp}),
\[\frac{1}{n}\sum W^2_t \stackrel{P_B}{\longrightarrow} 1 \] 
This proves condition (6).\\

\noindent
\textbf{Condition(7)}\hspace{.1in}$\lim_{k\rightarrow \infty}$
$\limsup_{n\rightarrow\infty} \sqrt{E(W^2_tI_{(|W_t|>k)})}\;=\; 0$

\begin{eqnarray*}
E(W^2_tI_{|W_t|>k}) &=& \frac{1}{\sigma^2_n}E[(w_t-1)^2I_{(|w_t-1|>k\sigma_n)}] \\
&\le& \frac{1}{\sigma^2_n}[E(w_t-1)^4]^{\frac{1}{2}}\;[P(|w_t-1|>k\sigma_n)]^{\frac{1}{2}} \\
&\le& \frac{1}{\sigma^2_n}(M^{\frac{1}{2}}_{n4})(\frac{\sigma^2_n}{k^2\sigma^2_n})^{\frac{1}{2}} \\
&=& \frac{1}{k}(\;\frac{M_{n4}}{\sigma^4_n}\;)^{\frac{1}{2}}
\end{eqnarray*}
where $M_{n4} = E(w_t-1)^4$.
Therefore \[\lim_{k\rightarrow \infty}
\limsup_{n\rightarrow\infty} \sqrt{E(W^2_tI_{(|W_t|>k)})}\]
\[\le \lim_{k\rightarrow \infty} \limsup_{n\rightarrow\infty}\frac{1}{\sqrt{k}}(\;\frac{M_{n4}}{\sigma^4_n}\;)^{\frac{1}{4}}
= 0\] as both $M_{n4}$ and $\sigma^4_n$ are bounded (follows from (\ref{mom4})). 
\end{description}
\medskip
\textbf{Case 2.}\hspace{0.05in}$\mathbf{(w_1,w_2,\ldots,w_n)\hspace{0.1in}iid\;
(1,\sigma^2)}$

\smallskip
\noindent
Again we need to establish (A4), that is, verify conditions 6) and 7) in Result(\ref{pw}).

\smallskip
\noindent Condition 6) follows from WLLN.

\smallskip
\noindent To verify condition 7), note that since distribution of $(w_1,w_2,\ldots,w_n)$
is independent of n,
\[\lim_{k\rightarrow \infty}
\limsup_{n\rightarrow\infty} \sqrt{E(W^2_tI_{(|W_t|>k)})}\]
\[=\lim_{k\rightarrow \infty}\sqrt{E(W^2_tI_{(|W_t|>k)})}\;=0\]
since $EW^2_t < \infty$.\\ \\

\noindent
\textbf{Remark 1.} Result 2 is true even when the process is nonstationary.
This follows from the fact that, given observations $\{X_t\}$ from the AR process, $X_t = \theta X_{t-1} + Z_t$, $|\theta|<1$; we can get a stationary solution of the above process, say \{$Y_t$\}, such that $\frac{1}{n}\sum X^a_t Z^b_{t+k}\stackrel{a.s.}{\longrightarrow}E(Y^a_tZ^b_{t+k})$.\\
As a consequence, Theorem 1) holds even without the assumption of stationarity, which is assumed throughout its proof.

\section{Bootstrap in Heteroscedastic AR(1) model} Now we introduce 
heteroscedasticity in the model (\ref{ar}), and study the Weighted Bootstrap estimator.
Consider the following model:

\begin{eqnarray}
&& X_t = \theta_0 X_{t-1} + Z_t; \; Z_t = \tau_t \epsilon_t  
 \hspace{0.1in}t=1,2,\ldots,n.\hspace{.1in}|\theta_0|<1
 \label{mh} \\
&& X_0 \sim F_0 \; with\; all\; moments\; finite.
\end{eqnarray}

\noindent where $\theta_0$, $\tau_t > 0$ are constants,$\epsilon_t \sim iid(0,1)$, and $\epsilon_t$ is independent of $\{ X_{t-k}, k \ge 1 \}$ for all $t$. \\ \\

\subsection{Estimation}
Based on observations $X_1, X_2, \ldots, X_n$ we discuss various methods for estimating $\theta$ in the model. Listed below are four types of estimators.\\
\smallskip

\noindent
\textbf{(a) Weighted Least Squares Estimator}
Assuming $\{\tau_t\}$ to be known, consider the following estimator for $\theta_0$:
\begin{eqnarray}
\hat\theta_1 &=& argmin_{\theta} \frac{1}{n}\sum_{t=2}^n \frac{1}{\tau^2_t}(X_t - \theta X_{t-1})^2 \\
&=& \frac{\sum_{t=2}^n \frac{1}{\tau^2_t} X_t X_{t-1}}{\sum_{t=2}^n \frac{1}{\tau^2_t} X^2_{t-1}} \label{eh1}
\end{eqnarray}
If $\epsilon_t$ in model(\ref{mh}) is normal, (\ref{eh1}) turns out to be the (Gaussian) maximum likelihood estimators.\\

\medskip
\noindent \textbf{(b) Least Squares Estimator}
In general $\{\tau_t\}$ are unknown and are non-estimable. Hence we may consider the general least squares estimators, ie,
\begin{equation}
\hat\theta_2 = \frac{\sum_{t=2}^n X_t X_{t-1}}{\sum_{t=2}^n X^2_{t-1}} \label{eh2}
\end{equation}
This turns out to be the same as (\ref{eh1}) if the $\{\tau_i\}$ are all equal, that is the model is homoscedastic.\\

\medskip
\noindent \textbf{(c) Weighted Least Absolute Deviations Estimator} The estimators (\ref{eh1}) and (\ref{eh2}) are $L_2$-estimators. It is well known that $L_1$-estimators are more robust with respect to heavy-tailed distributions than $L_2$-estimators. This motivates the study of various LAD estimators for $\theta_0$. Now we reparametrize model(\ref{mh}) in such a way that the median of $\epsilon_t$, instead of the mean equals 0 while $V\epsilon_t = 1$ remains unchanged. Our first absolute deviation estimator takes the form
\begin{equation}
\hat\theta_3 = argmin_{\theta} \sum_{t=2}^n \frac{1}{\tau_t} |X_t - \theta X_{t-1}| \label{eh3}
\end{equation}
This is motivated by the fact that $\hat\theta_3$ turns out to be the maximum likelihood estimator when the errors have double-exponential distribution.\\

\medskip
\noindent \textbf{Least absolute deviations estimator} Estimator(\ref{eh3}) uses the fact that $\tau_t$ are known. Incase they are not our absolute deviation estimator takes the form

\begin{equation}
\hat\theta_4 = argmin_{\theta} \sum_{t=2}^n |X_t - \theta X_{t-1}| \label{eh4}
\end{equation}
In the next section we discuss the assymptotic properties of the listed estimators.

\subsection{ Consistency of estimation in heteroscedastic AR(1) process}
In this section, we establish the distributional consistency of each of the four estimators discussed in the earlier section. To do so, we will use some established results, the first one being the following Martingale Central Limit theorem:

\begin{res}[Martingale C.L.T.; see Hall and Heyde 1980] \label{mclt}
 Let \{$S_n, \mathcal{F}_n$\} denote a zero-mean martingale whose increments have finite
 variance. Write $\;S_n = \sum_{i=1}^{n} X_i\;$, 
 $\;V^2_n = \sum_{i=1}^nE( X^2_{i-1}|\mathcal{F}_{i-1} )\;$
 and $\;s^2_n = EV^2_n = ES^2_n\;$. If
\begin{eqnarray*}
s^{-2}_n V^2_n &\stackrel{P}{\longrightarrow}& 1 \;\;and \\
s^{-2}_n\sum_{i=1}^{n}E(X^2_iI(|X_i|\ge \epsilon s_n)) &\longrightarrow& 0\;\; as \; n\rightarrow \infty \;\;\forall \epsilon > 0.
\end{eqnarray*}
 Then $\frac{S_n}{s_n}\stackrel{d}{\longrightarrow}N(0,1)$.
\end{res}

\NI
Another result we will need is the following one on convergence of a weighted sum of iid random variables.

\begin{res} \label{wcp}
Let $X_1, X_2, \ldots, X_n$ be a sequence of iid mean zero random variables, and $\{c_{in} | i=1,\ldots,n \}$ a triangular sequence of bounded constants. Then $ \frac{1}{n} \sum_{i=1}^n c_{in} X_i \cp 0$
\end{res}

\subsubsection{Distributional consistency of $\hat\theta_1$}

\begin{thm} \label{th1}
Define $\frac{s^2_n}{n} = \frac{1}{n}\sum_{t=2,n} \tau_t^{-2} EX^2_{t-1}$. Assume that

\bed
\item[A1.] $\frac{\tau_i}{\tau_j} \le M_2 \; \forall \; 1 \le i < j \le n$
\item[A2.] $\frac{1}{n} \sum_{1 \le i < j \le n} \Th^{2(j-i)}_0 (\frac{\tau_i}{\tau_j})^2 \ge M_1 >0$ 
\item[A3.] $\frac{1}{n} \sum_{1 \le i < j \le n} \Th^{2(j-i)}_0 (\frac{\tau_i}{\tau_j})^2 \lra \rho^2$
\ed 

\NI
Then under assumptions(A1-A2), $s_n(\hat\Th_1 - \Th_0)\cd N(0,1)$\\
Further if we assume (A3), we have, $\sqrt{n}(\hat\Th_1 -
\Th_0)\cd N(0,\Th^2_0 / \rho^2)$ as $n\rightarrow\infty$ 
\end{thm}

\NI
\textbf{Proof}
$\sqrt{n}(\hat\Th_1-\Th_0) = \sqrt{n}\frac{\sum \tau^{-2}_t X_{t-1}Z_t}{\sum \tau^{-2}_t X^2_{t-1}}$ \\ 

\smallskip
\NI
\tb{Step1.} $\frac{1}{\sqrt{n}} \sum \tau^{-2}_t X_{t-1}Z_t$ is assymptotically normal \\

\NI
Let $S_n = \sum_{t=2}^{n} \tau^{-2}_t X_{t-1}Z_{t}$.\\ 
Note that 
\beq 
X_t &=& \Th^t_0X_0 + \sum_{k=1}^t\Th^{t-k}_0Z_k\hspace{.1in}\forall t \; \ge 1 \label{x} \\
Hence \; E(X^2_t) &=& \Th^{2t}_0EX^2_0 + \sum_{k=1}^{t}\Th^{2t - 2k}_0\tau^2_k \label{ex2}
\eeq

\NI
Hence $S_n$ is a 0 mean $\mathcal{A}_n$ measurable martingale with increments having finite variance, where $\mathcal{A}_t = \sigma(X_0, \eps_1, \eps_2,\ldots,\eps_t);\hspace{0.1in} t=1,2,\ldots,n$. This follows from the fact that $E(X^2_t)$ is finite, and $E(X_{t-1}Z_t|\mathcal{A}_{t-1}) = 0$.\\
To establish the assymptotic normality of $S_n$, we use Result (\ref{mclt}).
Let
\beqn
V_n^2 &=& \sum_{t=2}^{n}E(\tau^{-4}_t X^2_{t-1}Z^2_t|\mathcal{A}_{t-1}) \\
&=& \sum_{t=2}^{n}\tau^{-2}_t X^2_{t-1} 
\eeqn
Then to accomplish Step1, we need to show
\beq
&& \frac{V^2_n}{s^2_n} \cp 1 \label{mch1} \\
&& \frac{1}{s^2_n} \sum_{t=2,n} E \left[ (\tau^{-2}_t X_{t-1}Z_t)^2 I(\tau^{-2}_t |X_{t-1}Z_t| \ge \eps s_n) \right] \lra 0 \label{mch2}
\eeq

\NI
To prove (\ref{mch1}), note that

\beqn
&& \frac{V^2_n}{n} - \frac{s^2_n}{n} \\
&=& \frac{1}{n} \sum_{t=1,n-1} \tau^{-2}_{t+1} [X^2_t - EX^2_t] \\
&=& \frac{1}{n} \sum_{t=1}^{n-1} \sum_{k=0}^t \frac{\tau^2_k}{\tau^2_{t+1}} \Th_0^{2(t-k)}(\eps^2_k - 1) + \frac{2}{n}\sum_{t=1,n-1} \sum_{0 \le i < j \le t} \frac{\tau_i \tau_j}{\tau^2_{t+1}} \Th^{2t-i-j}_0 \eps_i \eps_j \\
&=& \frac{1}{n}\sum_{k=0}^{n-1}(\eps^2_k - 1)\left(\sum_{t=k}^{n-1}\frac{\tau^2_k}{\tau^2_{t+1}} \Th_0^{2(t-k)}\right) +
\frac{2}{n}\sum_{0 \le i < j \le n-1}\eps_i \eps_j\left(\sum_{t=j}^{n-1}\frac{\tau_i \tau_j}{\tau^2_{t+1}}\Th^{2t-i-j}_0 \right) \\
&=& T_1 + 2T_2 \; (say)
\eeqn

\NI
Using assumption (A1) and Resut\ref{wcp}, we have $T_1 \cp 0$.

\beqn
&& ET_2 = 0 \\
&& VT_2 = \frac{1}{n^2}\sum_{0 \le i < j \le n}\sum_{t=j}^n \Th^{2t-i-j}_0 \frac{\tau_i \tau_j}{\tau^2_{t+1}} \\
&& \le \frac{M_2}{n^2} \sum_{t=1}^{n-1} \left(\sum_{k=0}^t |\Th_0|^{t-k}\right)^2 \lra 0
\eeqn
Hence $T_2 \cp 0$.\\
Combining, $\frac{V^2_n}{n} - \frac{s^2_n}{n} \cp 0$\\
Also $\frac{s^2_n}{n} = \frac{1}{n} \sum_{t=1}^{n-1} \sum_{k=0}^t \frac{\tau^2_k}{\tau^2_{t+1}}\Th^{2(t-k)}_0$.\\
Using assumption(A2), $\frac{s^2_n}{n}$ is bounded below.
Therefore $\frac{V^2_n}{s^2_n} \cp 1$ and this proves (\ref{mch1}).\\

\smallskip
\NI 
Remains to show (\ref{mch2}), ie $\frac{1}{s^2_n} \sum_{t=1}^{n-1} E \left[ \tau^{-2}_{t+1} X^2_t \eps_{t+1}^2 I(\tau^{-1}_{t+1} |X_t \eps_{t+1}| \ge \eps s_n) \right] \lra 0$ \\

\beqn
&& |X_t| \le \sum_{k=0}^t \tau_k |\Th_0|^{t-k} |\eps_k| \\
&& \frac{X^2_t}{\tau^2_{t+1}} \le \sum_{k=0}^t \frac{\tau^2_k}{\tau_{t+1}^2} |\Th_0|^{2(t-k)}\eps_k^2 + 2\sum_{0\le i < j \le t} \frac{\tau_i \tau_j}{\tau_{t+1}^2} |\Th_0|^{2t-i-j} \eps_i \eps_j \\
\eeqn
Hence for $1 \le t \le n$,
\beqn
&& E \left[ \tau^{-2}_{t+1} X^2_t \eps_{t+1}^2 I(\tau^{-1}_{t+1} |X_t \eps_{t+1}| \ge \eps s_n) \right] \\
&\le& \sum_{k=0}^t \frac{\tau_k^2}{\tau_{t+1}^2} |\Th_0|^{2(t-k)} E \left[ \eps_k^2 \eps_{t+1}^2 I(\tau^{-1}_{t+1} |X_t \eps_{t+1}| \ge \eps s_n) \right] + 2\sum_{0 \le i < j \le t}\frac{\tau_i \tau_j}{\tau_{t+1}^2} |\Th_0|^{2t-i-j} E\left[\eps_i\eps_j \eps_{t+1}^2 I(\tau^{-1}_{t+1} |X_t \eps_{t+1}| \ge \eps s_n) \right] \\
&\le& A_1 \max_{0 \le k \le t} E \left[ \eps_k^2 \eps_{t+1}^2
I(\tau^{-1}_{t+1} |X_t \eps_{t+1}| \ge \eps s_n) \right] 
+ A_2 \max_{0 \le i < j \le t} E\left[\eps_i\eps_j \eps_{t+1}^2 I(\tau^{-1}_{t+1} |X_t \eps_{t+1}| \ge \eps s_n) \right] \\
&\le& A \max_{0 \le k \le t} E \left[ \eps_k^2 \eps_{t+1}^2 I(\tau^{-1}_{t+1} |X_t \eps_{t+1}| \ge \eps s_n) \right]
\eeqn 

\NI $\frac{s_n}{\sn}$ is bounded below by say, $M > 0$. Hence for a fixed $k_0$, $0 \le k_0 \le t$,

\beqn
&& E \lt[ \eps_{k_0}^2 \eps_{t+1}^2 I \lt(\frac{|X_t \eps_{t+1}|}{\tau_{t+1}} \ge \eps s_n \rt) \rt] \\
&\le& E \left[ \eps_{k_0}^2 \eps_{t+1}^2 I\left(\frac{|X_t \eps_{t+1}|}{\tau_{t+1}} \ge \eps M \sn \right) \right] \\
&\le& E\left[ \eps_{k_0}^2 \eps_{t+1}^2 I\left(|\eps_{t+1}|   
\sum_{k=0}^t \frac{\tau_k}{\tau_{t+1}}|\eps_k| \ge \eps M \sn \rt) \rt] \\
&\le& E\left[ \eps_{k_0}^2 \eps_{t+1}^2 I\left(|\eps_{t+1}|   
\ge \sqrt{\eps M} n^{1/4} \rt) \rt] + E\left[ \eps_{k_0}^2 \eps_{t+1}^2 I\left(\sum_{k=0}^t \frac{\tau_k}{\tau_{t+1}}
|\Th_0|^{t-k}|\eps_k| \ge \sqrt{\eps M}n^{1/4} \rt) \rt] \\
&=& E\left[\eps_{t+1}^2 I\left(|\eps_{t+1}|\ge \sqrt{\eps M} n^{1/4} \rt) \rt] + E\left[ \eps_{k_0}^2 I\left(\sum_{k=0}^t \frac{\tau_k}{\tau_{t+1}}|\Th_0|^{t-k}|\eps_k| \ge \sqrt{\eps M}n^{1/4} \rt) \rt] \\
&\le& E\left[\eps_1^2 I\left(|\eps_1|\ge c_1 n^{1/4} \rt) \rt] + E\left[ \eps_{k_0}^2 I\left(\sum_{k=0}^t |\Th_0|^{t-k}|\eps_k| \ge c_2 n^{1/4} \rt) \rt] \\
&\le& E\left[\eps_1^2 I\left(|\eps_1|\ge c_1 n^{1/4} \rt) \rt] + 
E\lt[ \eps_{k_0}^2 I\left(|\eps_{k_0}| \ge \frac{c_2}{2} n^{1/4} \rt) \rt] +  
E\lt[ \eps_{k_0}^2 I\left(\sum_{k\ne k_0}|\Th_0|^{t-k}|\eps_k| \ge \frac{c_2}{2} n^{1/4} \rt) \rt] \\
&\le& E\left[\eps_1^2 I\left(|\eps_1|\ge c_1 n^{1/4} \rt) \rt] + 
E\lt[ \eps_1^2 I\left(|\eps_1| \ge \frac{c_2}{2} n^{1/4} \rt) \rt] +  
P\lt[\sum_{k\ne k_0}|\Th_0|^{t-k}|\eps_k| \ge \frac{c_2}{2} n^{1/4} \rt] \\
&\le& E\left[\eps_1^2 I\left(|\eps_1|\ge c_3 n^{1/4} \rt) \rt]
+ \frac{c_4}{n^{1/4}} 
\eeqn
Hence $\max_{0 \le k \le t} E \left[ \eps_k^2 \eps_{t+1}^2 I(\tau^{-1}_{t+1} |X_t \eps_{t+1}| \ge \eps s_n) \right] \lra 0$. Using the fact that $\frac{s_n^2}{n}$ is bounded below, this proves (\ref{mch2}).
Using Result\ref{mclt}, from (\ref{mch1}) and (\ref{mch2}) we deduce that $\frac{S_n}{s_n} \cd N(0,1)$, ie
\be
\frac{1}{s_n} \sum_{t=2}^n \tau_t^{-2}X_{t-1}Z_t \cd N(0,1) \label{s1}
\ee

\NI
\tb{Step2.} 
\be
\frac{\sum_{t=2}^n \tau_t^{-2}X_{t-1}^2}{s_n^2} \cp 1 \label{s2}
\ee
This follows from (\ref{mch1}).\\

\smallskip
\NI
From (\ref{s1}) and (\ref{s2}) we deduce,

\beq 
\NI s_n \frac{ \sum \taur X_{t-1}Z_t}{\sum \taur X_{t-1}^2} \cd N(0,1) \\
\NI ie\; s_n(\hat\Th_1 - \Th_0) \cd N(0,1)
\eeq

\beq
\frac{s_n^2}{n} &=& \frac{1}{n} \sum_{t=2}^n \taur EX_{t-1}^2 \\
&=& \frac{1}{n} \sum_{t=1}^{n-1} \sum_{k=1}^t \frac{\tau_k^2}{\tau_{t+1}^2}\Th_0^{2(t-k)}
\eeq
Hence if we assume (A3), we have $\frac{s_n^2}{n} \ra \frac{\rho^2}{\Th_0^2}$, and then,
\[
\mb{ \sn(\hat\Th_1 - \Th_0) \cd N(0, \frac{\Th_0^2}{\rho^2}) }
\]
This completes the proof. \\

\smallskip
\NI
\tb{Remark} Assumptions (A1) and (A2) are satisfied if     
$\{\tau_t\}$'s are bounded, or more generally if they are of the same order, ie there exists constants $c_1, c_2 > 0$ and 
$\alpha \ge 0$ such that $c_1t^{\alpha} \le \tautsq \le c_2t^{\alpha}$ for $1\le t \le n$.

\subsubsection{Distributional consistency of $\hat\theta_2$}

\begin{thm} \label{th2}
Define $s_n^2 = \sum_{t=2}^n\tautsq E(X_{t-1}^2)$. Suppose
$\{\tau_t\}$'s satisfy the following assumptions.
\bed
\item[A1.] $M_1 \le \tau_t \le M_2;\; t=1,2,\ldots,n$
\item[A2.] $\frac{\sum\tautsq}{n} \ra \tausq > 0$
\item[A3.] $\frac{1}{n}\sum_{1\le i<j \le n}\tausq_i  
\tausq_j\Th_0^{2(j-i)} \ra \rho^2$
\ed 
Then under assumptions(A1) and (A2), $\lt(\frac{s_n}{\sn}\rt)^{-1}\sn(\hat\Th_2 - \Th_0) \cd  N\lt(0 , \frac{(1-\Th_0^2)^2}{\tau^4}\rt)$.\\
Further if (A3) holds, $\sn(\hat\Th_2 - \Th_0) \cd
N\lt(0 , \frac{\rho^2(1-\Th_0^2)^2}{\tau^4\Th_0^2}\rt)$.
\end{thm}

\NI
\textbf{Proof}
$\sqrt{n}(\hat\Th_2-\Th_0) = \sqrt{n}\frac{\sum X_{t-1}Z_t}{\sum X^2_{t-1}}$ \\

\smallskip
\NI
\tb{Step1.} $\frac{1}{\sqrt{n}}\sum X_{t-1}Z_t$ is assymptotically normal. \\

\NI
Let $S_n = \sum_{t=2}^{n}X_{t-1}Z_{t}$.\\ 

\NI
Then $S_n$ is a 0 mean $\mathcal{A}_n$ measurable martingale with increments having finite variance, where $\mathcal{A}_t = \sigma(X_0, \eps_1, \eps_2,\ldots,\eps_t);\hspace{0.1in} t=1,2,\ldots,n$. This follows from the fact that $E(X^2_t)$ is finite, and $E(X_{t-1}Z_t|\mathcal{A}_{t-1}) = 0$.\\
To establish the assymptotic normality of $S_n$, we use Result (\ref{mclt}).
Let
\beqn
V_n^2 &=& \sum_{t=2}^{n}E(X^2_{t-1}Z^2_t|\mathcal{A}_{t-1}) \\
&=& \sum_{t=2}^{n}\tautsq X^2_{t-1} 
\eeqn
Then to accomplish Step1, we need to show
\beq
&& \frac{V^2_n}{s^2_n} \cp 1 \label{mc21} \\
&& \frac{1}{s^2_n} \sum_{t=2}^n E\lt[(X_{t-1}Z_t)^2 
I(|X_{t-1}Z_t|\ge\eps s_n) \rt] \ra 0 \label{mc22}
\eeq

\NI
Using the expressions for $X_t$ and $EX_t^2$ from equations(\ref{x}) and (\ref{ex2}), we have

\beqn
&& \frac{V^2_n}{n} - \frac{s^2_n}{n} \\
&=& \frac{1}{n}\sum_{t=1}^{n-1}\tausq_{t+1}[X^2_t - EX^2_t] \\
&=& \frac{1}{n}\sum_{t=1}^{n-1}\sum_{k=0}^t \tausq_{t+1} 
\tausq_k \Th_0^{2(t-k)}(\eps^2_k - 1) + \frac{2}{n}\sum_{t=1}^{n-1}\sum_{0 \le i < j \le t}\tau_i\tau_j\tau^2_{t+1}
\Th^{2t-i-j}_0 \eps_i \eps_j \\
&=& \frac{1}{n}\sum_{k=0}^{n-1}(\eps^2_k - 1)\lt(\sum_{t=k}^{n-1}\tau^2_k\tau^2_{t+1}\Th_0^{2(t-k)}\rt) +
\frac{2}{n}\sum_{0\le i<j\le n-1}\eps_i\eps_j\lt(\sum_{t=j}^{n-1}\tau_i\tau_j\tau^2_{t+1}\Th^{2t-i-j}_0 \rt) \\
&=& T_1 + 2T_2 \; (say)
\eeqn

\NI
Using assumption (A1) and Resut\ref{wcp}, we have $T_1 \cp 0$.
\beq
&& ET_2 = 0 \\
&& VT_2 = \frac{1}{n^2}\sum_{0 \le i < j \le n-1}\lt(\sum_{t=j}^{n-1} \Th^{2t-i-j}_0\tau_i\tau_j\tau^2_{t+1}\rt)^2 \label{l2} \\
&\le& \frac{c}{n^2}\sum_{0\le i<j\le n-1}\sum_{t=j}^{n-1}
|\Th_0|^{2t-i-j} \label{l3} \\
&\le& \frac{c}{n^2} \sum_{t=1}^{n-1} \left(\sum_{k=0}^t |\Th_0|^{t-k}\right)^2 \lra 0
\eeq
Here $c$ is some positive constant. (\ref{l3}) follows from (\ref{l2}) using the fact that $\sum_{t=j}^{n-1} \Th^{2t-i-j}_0\tau_i\tau_j\tau^2_{t+1}$ is bounded which inturn follows from assumption(A1).\\ 
Hence $T_2 \cp 0$.\\
Combining, $\frac{V^2_n}{n} - \frac{s^2_n}{n} \cp 0$\\
Also $\frac{s^2_n}{n} = \frac{1}{n} \sum_{t=1}^{n-1} \sum_{k=0}^t \frac{\tau^2_k}{\tau^2_{t+1}}\Th^{2(t-k)}_0$.\\
Again using assumption(A1), $\frac{s^2_n}{n}$ is bounded below.
Therefore $\frac{V^2_n}{s^2_n} \cp 1$ and this proves (\ref{mc21}).\\

\smallskip
\NI 
Remains to show (\ref{mc22}), ie $\frac{1}{s^2_n} \sum_{t=1}^{n-1} E \left[ \tau^{2}_{t+1} X^2_t \eps_{t+1}^2 I(\tau_{t+1} |X_t \eps_{t+1}| \ge \eps s_n) \right] \lra 0$ \\
\beqn
&& |X_t| \le \sum_{k=0}^t \tau_k |\Th_0|^{t-k} |\eps_k| \\
&& X^2_t \le \sum_{k=0}^t\tau^2_k|\Th_0|^{2(t-k)}\eps_k^2 + 2\sum_{0\le i<j \le t}\tau_i\tau_j|\Th_0|^{2t-i-j}\eps_i\eps_j \\
\eeqn
Hence for $1 \le t \le n$,
\beqn
&& E \lt[ \tau^2_{t+1} X^2_t \eps_{t+1}^2 I(\tau_{t+1} |X_t \eps_{t+1}| \ge \eps s_n) \rt] \\
&\le& \sum_{k=0}^t\tau_k^2\tau_{t+1}^2|\Th_0|^{2(t-k)}E\lt[\eps_k^2\eps_{t+1}^2 I(\tau_{t+1}|X_t \eps_{t+1}| \ge \eps s_n) \rt] + 2\sum_{0 \le i < j \le t}\tau_i \tau_j\tau_{t+1}^2 |\Th_0|^{2t-i-j} E\lt[\eps_i\eps_j \eps_{t+1}^2 I(\tau_{t+1} |X_t \eps_{t+1}| \ge \eps s_n) \rt] \\
&\le& A_1 \max_{0 \le k \le t} E \lt[\eps_k^2 \eps_{t+1}^2
I(\tau_{t+1} |X_t \eps_{t+1}| \ge \eps s_n) \rt] 
+ A_2 \max_{0 \le i < j \le t} E\left[\eps_i\eps_j \eps_{t+1}^2 I(\tau_{t+1} |X_t \eps_{t+1}| \ge \eps s_n) \right] \\
&\le& A \max_{0 \le k \le t} E \left[ \eps_k^2 \eps_{t+1}^2 I(\tau_{t+1} |X_t \eps_{t+1}| \ge \eps s_n) \rt]
\eeqn 

\NI $\frac{s_n}{\sn}$ is bounded below by say, $M > 0$. Hence for a fixed $k_0$, $0 \le k_0 \le t$,

\beqn
&& E \lt[ \eps_{k_0}^2 \eps_{t+1}^2 I \lt(\tau_{t+1}|X_t \eps_{t+1}|\ge \eps s_n \rt) \rt] \\
&\le& E \left[ \eps_{k_0}^2 \eps_{t+1}^2 I\lt(\tau_{t+1}
|X_t \eps_{t+1}|\ge \eps M \sn \rt) \rt] \\
&\le& E\lt[ \eps_{k_0}^2 \eps_{t+1}^2 I\left(|\eps_{t+1}|   
\sum_{k=0}^t\tau_k\tau_{t+1}|\eps_k| \ge \eps M \sn \rt) \rt] \\
&\le& E\left[ \eps_{k_0}^2 \eps_{t+1}^2 I\left(|\eps_{t+1}|   
\ge \sqrt{\eps M} n^{1/4} \rt) \rt] + E\left[ \eps_{k_0}^2 \eps_{t+1}^2 I\left(\sum_{k=0}^t\tau_k\tau_{t+1}
|\Th_0|^{t-k}|\eps_k| \ge \sqrt{\eps M}n^{1/4} \rt) \rt] \\
&=& E\left[\eps_{t+1}^2 I\left(|\eps_{t+1}|\ge \sqrt{\eps M} n^{1/4} \rt) \rt] + E\left[ \eps_{k_0}^2 I\left(\sum_{k=0}^t \tau_k\tau_{t+1}|\Th_0|^{t-k}|\eps_k| \ge \sqrt{\eps M}n^{1/4} \rt) \rt] \\
&\le& E\lt[\eps_1^2 I\left(|\eps_1|\ge c_1 n^{1/4} \rt) \rt] + E\lt[ \eps_{k_0}^2 I\left(\sum_{k=0}^t |\Th_0|^{t-k}|\eps_k| \ge c_2 n^{1/4} \rt) \rt] \\
&\le& E\left[\eps_1^2 I\left(|\eps_1|\ge c_1 n^{1/4} \rt) \rt] + 
E\lt[ \eps_{k_0}^2 I\left(|\eps_{k_0}| \ge \frac{c_2}{2} n^{1/4} \rt) \rt] +  
E\lt[ \eps_{k_0}^2 I\left(\sum_{k\ne k_0}|\Th_0|^{t-k}|\eps_k| \ge \frac{c_2}{2} n^{1/4} \rt) \rt] \\
&\le& E\left[\eps_1^2 I\left(|\eps_1|\ge c_1 n^{1/4} \rt) \rt] + 
E\lt[ \eps_1^2 I\left(|\eps_1| \ge \frac{c_2}{2} n^{1/4} \rt) \rt] +  
P\lt[\sum_{k\ne k_0}|\Th_0|^{t-k}|\eps_k| \ge \frac{c_2}{2} n^{1/4} \rt] \\
&\le& E\left[\eps_1^2 I\left(|\eps_1|\ge c_3 n^{1/4} \rt) \rt]
+ \frac{c_4}{n^{1/4}} 
\eeqn
Hence $\max_{0 \le k \le t} E \left[ \eps_k^2 \eps_{t+1}^2 I(\tau_{t+1} |X_t \eps_{t+1}| \ge \eps s_n) \right] \lra 0$. Using the fact that $\frac{s_n^2}{n}$ is bounded below, this proves (\ref{mc22}).\\
Using Result\ref{mclt}, from (\ref{mc21}) and (\ref{mc22}) we deduce that $\frac{S_n}{s_n} \cd N(0,1)$, ie
\be
\frac{1}{s_n} \sum_{t=2}^n X_{t-1}Z_t \cd N(0,1) \label{s21}
\ee

\NI
\tb{Step2.} 
\be
\IN\sum_{t=2}^n X_{t-1}^2 \cp \frac{\tausq}{(1-\Th_0^2)}
\label{s22}
\ee
This follows once we show
\beq
&& \IN\sum_{t=1}^{n-1}(X_t^2 - EX_t^2) \cp 0 \label{st21} \\
and \; &&\IN\sum_{t=1}^{n-1}EX_t^2 \ra \frac{\tausq}{(1-\Th_0^2)} \label{st22}
\eeq
Using the expressions for $X_t$ and $EX_t^2$ from equations(\ref{x}) and (\ref{ex2}), we have
\beqn
&& \IN\sum_{t=1}^{n-1}(X_t^2 - EX_t^2) \\ 
&=& \IN\sum_{k=0}^{n-1}(\eps^2_k - 1)\lt(\sum_{t=k}^{n-1}\tau^2_k\Th_0^{2(t-k)}\rt) + \frac{2}{n}\sum_{0\le i<j\le n-1}\eps_i\eps_j\lt(\sum_{t=j}^{n-1}\tau_i\tau_j\Th^{2t-i-j}_0 \rt) \\
&\cp 0
\eeqn
The above steps can be justified by proceeding as in the proof of (\ref{mc21}). This complete (\ref{st21}). \\
To see (\ref{st22}), note that
\beq
\IN\sum_{t=1}^{n-1}EX_t^2 &=& \frac{EX_0^2}{n}\sum_{t=1}^{n-1}\Th_0^{2t} + \IN\sum_{t=1}^{n-1}\sum_{k=1}^t\tausq_k\Th_0^{2t-2k} \\
&\approx& \IN\sum_{t=1}^{n-1}\sum_{k=1}^t\tausq_k\Th_0^{2t-2k} \\
&=& \IN\sum_{k=1}^{n-1}\tausq_k\sum_{t=0}^{n-k-1}\Th_0^{2t} \\
&=& \IN\sum_{k=1}^{n-1}\tausq_k\frac{(1-\Th_0^{2(n-k)})}{(1-\Th_0^2)} \\
&\approx& \frac{1}{(1-\Th_0^2)}\IN\sum\tausq_k \\
&\ra& \frac{\tausq}{(1-\Th_0^2)}
\eeq
This proves (\ref{st22}). (\ref{st21}) and (\ref{st22}) together prove (\ref{s22}) and this completes Step2.\\

\smallskip
\NI
Dividing (\ref{s21}) by (\ref{s22}) we deduce,
\[
(\frac{s_n}{\sn})^{-1}\sn(\hat\Th_2 - \Th_0) \cd N\lt(0,\frac{(1-\Th_0^2)^2}{\tau^4}\rt)
\]

\beqn
\frac{s_n^2}{n} &=& \IN\sum_{t=2}^n\tau_t^2 EX_{t-1}^2 \\
&=& \IN\sum_{t=1}^{n-1} \sum_{k=1}^t\tau_k^2\tau_{t+1}^2
\Th_0^{2(t-k)}
\eeqn
Hence if we assume (A3), we have $\frac{s_n^2}{n} \ra \frac{\rho^2}{\Th_0^2}$, and then,
\[
\mb{ \sn(\hat\Th_2 - \Th_0) \cd N\lt(0, \frac{\rho^2(1-\Th_0^2)^2}{\tau^4\Th_0^2}\rt) }
\]
This completes the proof. \\

\subsection{ Consistency of the Weighted Bootstrap technique}
Now we prove that the Weighted Bootstrap estimator is assymptotically normal with the same assymptotic
distribution as of the least squares estimate. In particular WB provides a consistent resampling
scheme in the AR model with introduced heteroscedasticity.

\noindent
\begin{thm} \label{tb2}
Let $\hat\theta^{*}_n$ be the weighted bootstrap estimator of $\hat\theta_n$ as defined in (\ref{wb1}). Then
under the conditions (A1)-(A4) on the weights,

\begin{equation}
P\left [\;\sqrt{n}\sigma^{-1}_n(\hat\theta^{*}_n - \hat\theta_n) \le x\;
\vline X_1,\ldots,X_n\right ] - P\left [Y \le x\right ] =
 o_P(1)\hspace{.1in} \forall \; x \in \mathbb{R}
\end{equation}
where $Y \sim N(0, \sigma^2)$, $\sigma^2$ being defined in Theorem (\ref{arh}).
\end{thm}

\noindent
\textbf{Proof} As in (\ref{expand}),
\begin{eqnarray}
 & & \sqrt{n}\sigma^{-1}_n (\hat\theta^*_n-\hat{\theta}_n) \nonumber \\
 &=& \frac{\sum W_t X_{t-1} Z_t/\sqrt{n}}{\sum w_t X^2_{t-1}/n} -
   \sqrt{n}(\hat{\theta}_n - \theta)\sigma^{-1}_n \frac{\sum (w_t-1) 
	X^2_{t-1}/n}{\sum w_t X^2_{t-1}/n} \nonumber \\
&=& T_1 - T_2 \label{split}
\end{eqnarray}

\noindent
\textbf{Claim 1.} There exists $\tau > 0$ such that
\[
P\left [\;\frac{1}{\sqrt{n}}\sum_{t=2}^n W_t X_{t-1}Z_t \le x\;
\vline X_1,\ldots, X_n \right ]
\stackrel{P}{\longrightarrow}\Phi(\;\;\frac{x}{\tau}\;\;) \hspace{.1in}
\forall\; x \in \mathbb{R}
\]
To see this let us verify the first three conditions of Result(\ref{pw})
with $c_{nj} = X_jZ_{j+1}$ and $U_{nj} = W_j$ for $j = 1,\ldots,{n-1}$.

\smallskip
\noindent
\textbf{Condition 1}\hspace{.1in}
$\frac{1}{n} \sum_{t=2}^n X_{t-1}Z_t\stackrel{P}{\longrightarrow}0$

\smallskip
\noindent Follows from Theorem \ref{arh}.\\ \\
\textbf{Condition 2}\hspace{.1in} $\frac{1}{n}\sum_{t=2}^n (X_{t-1}Z_t)^2\stackrel{P}{\longrightarrow}\tau^2$\\
\begin{eqnarray*}
Let\hspace{.1in}S_n &=& \frac{1}{n}\sum_{t=2}^{n}X^2_{t-1}Z^2_t \\
U_n &=& \frac{1}{n}\sum_{t=2}^{n}X^2_{t-1}(Z^2_t - \sigma^2_t) \\
V_n &=& \frac{1}{n}\sum_{t=2}^{n}X^2_{t-1}\sigma^2_t 
\end{eqnarray*}
Then $S_n = U_n + V_n$.

\begin{eqnarray*}
E(U_n)  &=& 0 \\
V(nU_n) &=& \sum_{t=2}^{n}V(X^2_{t-1}(Z^2_t - \sigma^2_t)) +
2\sum_{2\le s\;}\sum_{<\;t\le n}
Cov(X^2_{t-1}(Z^2_t - \sigma^2_t),X^2_{s-1}(Z^2_s - \sigma^2_s)) \\
&=& \sum_{t=2}^{n}V(X^2_{t-1}(Z^2_t - \sigma^2_t))\\
&\le& M\sum_{t=2}^{n}EX^4_{t-1}\;\; where \; EZ^4_t\le M
\end{eqnarray*}

\begin{eqnarray*}
Therefore\hspace{.1in}V(U_n)
&\le& \frac{M}{n^2}\sum EX^4_{t-1} \longrightarrow0\hspace{.1in} since\hspace{.05in} \frac{1}{n}\sum EX^4_{t-1} \hspace{.05in}is\hspace{.05in}bounded.
\end{eqnarray*}
Hence $U_n\stackrel{P}{\longrightarrow}0$. \label{tbr1}
\begin{eqnarray*}
V_n &=& \frac{1}{n}[\;\sigma^2_1\sum_{t \; odd}X^2_{t-1} +
\sigma^2_2 \sum_ {t \; even}X^2_{t-1} ]\\
&=& \frac{1}{2}\left [\; \sigma^2_1\frac{\sum_{t\; odd}X^2_{t-1}}{n/2}\; +\;
\sigma^2_2\frac{\sum_{t\; even}X^2_{t-1}}{n/2}\right ]\\
&\stackrel{P}{\longrightarrow}& 
\frac{1}{2}[\; \sigma^2_1\sigma^2_{e}(1 - \theta^4)^{-1} + 
\sigma^2_2\sigma^2_{o}(1 - \theta^4)^{-1}\;]\\
&=& \frac{1}{2(1-\theta^4)}[ \sigma^2_1\sigma^2_{e} + \sigma^2_2\sigma^2_{o}] \label{tbr2}
\end{eqnarray*}
Therefore
\begin{eqnarray*}
S_n &\stackrel{P}{\longrightarrow}& \frac{1}{2}\frac{[\sigma^2_1\sigma^2_{e} + \sigma^2_2\sigma^2_{o}]}{(1-\theta^4)}\\
&=& \frac{\theta^2(\sigma^4_1 + \sigma^4_2) + 2\sigma^2_1\sigma^2_2}{2(1-\theta^4)}\; = \; \tau^2
\end{eqnarray*}

\noindent
\textbf{Condition 3}\hspace{.1in} 
$n^{-1}\max(X^2_{t-1}Z^2_t)\stackrel{P}{\longrightarrow}0$

\smallskip
\noindent Given $\epsilon$ positive,
\begin{eqnarray*}
P(n^{-1}\max(X^2_{t-1}Z^2_t)>\epsilon) &\le& \sum_{2}^{n}P(Y_t > n\epsilon)\hspace{.1in}where\hspace{.05in}Y_t = X^2_{t-1}Z^2_t \\
&\le& \frac{1}{n^2\epsilon^2}\sum_{t=2}^{n}E(X^4_{t-1}Z^4_t) \\
&=& \frac{M}{n^2\epsilon^2}\sum_{t=2}^{n}EX^4_{t-1}\\
&\longrightarrow&0\hspace{.1in}since\hspace{.1in}\frac{1}{n}\sum EX^4_t\hspace{.1in}is\hspace{.1in}bounded.
\end{eqnarray*}
This proves Claim 1.\\  
\noindent
\[Hence\;\;for\;\;\tau^2 =\frac{\theta^2(\sigma^4_1 + \sigma^4_2) + 
2\sigma^2_1\sigma^2_2}{2(1-\theta^4)}\;,\]

\begin{equation} \label{tbr3}
\mathbf{
P_B\left[\;\frac{1}{\sqrt{n}}\sum_{t=2}^n W_t X_{t-1}Z_t \le x \right]
\stackrel{P}{\longrightarrow}\Phi(\;\;\frac{x}{\tau}\;\;) \hspace{.1in}
\forall \; x \in \mathbb{R}
}
\end{equation}

\noindent
\textbf{Claim 2.}
\[With\;\; c = \frac{1}{2}\frac{(\sigma^2_1 + \sigma^2_2)}{(1-\theta^2)}\]

\begin{equation} \label{tbr4}
P_B\left [\;\;\vline\;\frac{1}{n}\sum_{t=2}^{n} w_t X^2_{t-1} - c\;\vline
\; > \;\epsilon\;\;\right ]
\stackrel{P}{\longrightarrow}0\hspace{.1in}\forall \; \epsilon > 0
\end{equation}
Note that
\[
\frac{1}{n}\sum w_tX^2_{t-1} = 
\frac{1}{2}\left [\frac{\sum_{t\; odd}w_tX^2_{t-1}}{n/2} +
\frac{\sum_{t\;even}w_tX^2_{t-1}}{n/2} \right ]
\]
Using the fact that $\{X_t\}_{t\;even}$ and $\{X_t\}_{t\;odd}$ form two
homoscedastic AR(1) processes, from Claim 2(Theorem \ref{t1}) and Remark 1, we
get,
\[\frac{\sum_{t\; odd}w_tX^2_{t-1}}{n/2}\stackrel{P_B}{\longrightarrow}
\frac{\theta^2\sigma^2_1 + \sigma^2_2}{1-\theta^4}\hspace{.1in}a.s.\]and
\[\frac{\sum_{t\;even}w_tX^2_{t-1}}{n/2}\stackrel{P_B}{\longrightarrow}
\frac{\sigma^2_1 + \theta^2\sigma^2_2}{1-\theta^4}\hspace{.05in}a.s.\]
Hence
\[\frac{1}{n}\sum w_tX^2_{t-1}\stackrel{P_B}{\longrightarrow}\frac{1}{2}
\frac{(\sigma^2_1 + \sigma^2_2)}{(1-\theta^2)}\hspace{.1in}a.s.\]
This proves Claim 2.\\ \\

\noindent
\textbf{Claim 3.}
\begin{equation} \label{tbr5}
P_B\left [\;\;\vline\;\frac{1}{n}\sum_{t=2}^{n} (w_t-1)X^2_{t-1}\;\vline\;
 > \;\epsilon\;\;\right ]
\stackrel{P}{\longrightarrow}0\hspace{.1in}\forall \; \epsilon > 0
\end{equation}
This follows from equations (\ref{tbr0}) and (\ref{tbr4}).\\ \\
Note that as defined in (\ref{split}), 
\[
\sqrt{n}\sigma^{-1}_n(\hat\theta^{*}_n\;-\;\hat\theta_n)\;=\;T_1\;-\;T_2
\]
Then from (\ref{tbr3}) and (\ref{tbr4}),
\[
P_B(T_1 \le x)\; - \; P(T \le x)\;=\; o_P(1)
\]
\[
where\;\; T\; \sim \; \left [\;\frac{1}{2}\frac{(\sigma^2_1 + \sigma^2_2)}
{(1 - \theta^2)}\;\right ]^{-1}\;
N\left (0, \frac{\sigma^2_1\sigma^2_2 + \theta^2(\sigma^4_1 +
\sigma^4_2)/2}{(1 - \theta^4)}\right )
\]
Moreover using equations (\ref{tbr4}) and (\ref{tbr5}), from Claim 3(Theorem
\ref{t1}), we get,
\[
P_B(\;|T_2 > \epsilon|\;)\;=\; o_P(1)\;\;\forall\;\epsilon>0
\] 
Combining
\begin{equation} \label{final}
P_B[\;\sqrt{n}\sigma^{-1}_n(\hat\theta^{*}_n - \hat\theta_n) \le x\;] - 
P[\;Y \le x\;] = o_P(1)\hspace{.1in} \forall x \in \mathbb{R}
\end{equation}
where 
\[
Y\;\sim \; N\left (\;0,\; 4\frac{(1 - \theta^2)}{(1 + \theta^2)}
\frac{\sigma^2_1\sigma^2_2 + 
\theta^2(\sigma^4_1 + \sigma^4_2)/2}{(\sigma^2_1 + \sigma^2_2)^2}\;\right )
\]
 and this completes the proof.\\ \\
 
\noindent \textbf{Remark 2.}
 In Theorems 1 and 3, we have established the consistency of the Weighted
 Bootstrap estimator in probability, ie we have proved,
 $\forall\;x\;\in\;\mathbb{R}$,
 \[
 P_B(\;\sqrt{n}\sigma^{-1}_n(\hat\theta^{*}_n - \hat\theta_n)\;\le\;x\;)\;-\;
 P(\;\sqrt{n}(\hat\theta_n - \theta)\;\le\;x\;)\;=\;o_P(1)
 \] 
 The same results can be achieved almost surely. One can prove that,
 $\forall\;x\;\in\;\mathbb{R}$,
 \[
 P_B(\;\sqrt{n}\sigma^{-1}_n(\hat\theta^{*}_n - \hat\theta_n)\;\le\;x\;)\;-\;
 P(\;\sqrt{n}(\hat\theta_n - \theta)\;\le\;x\;)
 \longrightarrow0\; a.s.
 \] 
 To prove this, one needs to verify the conditions of Result(1) almost surely,
 and replace all convergence of sample moments of \{$X_t$\} in probability,
 by almost sure convergence in the proofs.\\
 
\vspace{0.1in}

\section{Numerical Calculations}
In this section, we compare numerically the performance of the Weighted Bootstrap and Residual Bootstrap techniques for an heteroscedastic AR(1) model, and exhibit numerically, the consistency of the Weighted Bootstrap estimator.
We simulated 50 observations from the AR process, 
\[X_t = \theta X_{t-1} + Z_t,\; t = 1,2,\ldots,n.\]
where ${Z_t}$ is a sequence of independent Normal mean-zero random variables with $EZ^2_t = \sigma^2_1$ if $t$ is odd and $EZ^2_t = \sigma^2_2$ if $t$ is even. For simulation purpose, we used $\theta = 0.5$, $\sigma^2_1 = 1$, and $\sigma^2_2 = 2$.\\
The unknown $\theta$ is estimated by its LSE $\hat\theta_n$ which came to be $0.4418$.\\
Let $V_n = \sqrt{n}(\hat\theta_n - \theta)$ be the quantity of interest which is to be estimated using resampling techniques. Let $V^*_n = \sqrt{n}(\hat\theta^*_n - \hat\theta_n)$ denote its bootstrap estimate for two different bootstrap techniques: the Residual Bootstrap (which tacitly assumes that all the $Z_t$'s have same variance) and the Weighted Bootstrap. In case of WB, we used i.i.d Normal(1,1) weights. We used 200 simulations to estimate the distribution of $V^*_n$ in both the cases.\\
We performed the KS test to compare the distributions of $V_n$ and $V^*_n$. To estimate the distribution of $V_n$, we used 200 simulations from the above process. The results of the test are as follows.
\begin{center}
{\textbf{Two-Sample Kolmogorov-Smirnov Test}}\\
Data: $V_n$ and $V^*_n$\\
Alternative hypothesis:\\
cdf of $V_n$ does not equal the cdf of $V^*_n$ for at least one sample point\vspace{0.1in}
\begin{tabular}{|c||c||c|}
\hline
BS Technique & KS value & p-value \\
\hline
RB & 0.12 & 0.0945\\
WB & 0.1 & 0.234\\
\hline
\end{tabular}
\end{center}
Figure 1a) presents the estimated densities of $V_n$ and $V^*_n$, with
$\hat\theta^{*}_n$ being the residual bootstrap estimator, while Figure 1b)
presents the estimated densities with $\hat\theta^{*}_n$ being the weighted 
bootstrap estimator.
From the table it can be seen that both the estimators pass the test, 
but WB does reasonably better. This is also obvious from the density plots.\\ \\
Next we introduced more heteroscedasticity in the model. 
This time we took 
$\sigma^2_1$ to be $1$, and $\sigma^2_2$ as $10$. $\hat\theta_n$ came to
be $0.47083$. Again we estimate $V_n$ by $V^*_n$ and performed a KS test to
determine the goodness of the fit. Now the results are as follows:
\begin{center}
{\textbf{Two-Sample Kolmogorov-Smirnov Test}}\\
Data: $V_n$ and $V^*_n$\\
Alternative hypothesis:\\
cdf of $V_n$ does not equal the cdf of $V^*_n$ for at least one sample point\vspace{0.1in}
\begin{tabular}{|c||c||c|}
\hline
BS Technique & KS value & p-value \\
\hline
RB & 0.135 & 0.0431\\
WB & 0.125 & 0.0734\\
\hline
\end{tabular}
\end{center}
Figure 2a) presents the estimated densities of $V_n$ and $V^*_n$ for RB, while Figure 2b) presents the estimated densities for WB.
From the table, it can be seen that RB fails. This is expected since it is not adapted for heteroscedasticity. It fails to capture the true model in such a situation. WB still performs well, but its performance also falls. This is also reflected frm the density plots. Perhaps a larger sample size is required in case of substantial heteroscedasticity.\\
This illustrates the point that for small sample sizes, 
at small levels of heteroscedasticity, many Bootstrap techniques perform well
, but at substantial levels a careful choice is needed. The success of WB for
both levels of heteroscedasticity lends further support to our theoretical results.

\section{ARCH models} In this section, we first present the basic
probabilistic properties of ARCH models. Then we introduce various estimation
procedures for the parameters involved, and study their properties. The
assymptotic properties of the listed estimators under different error
distributions are also introduced. To approximate the distribution of the
estimators and draw inference based on an observed sample, various resampling
techniques are also listed along with their properties. Finally we supplement
our theoretical results with numerical calculations based on a simulated ARCH
data set.
 
\subsection{Basic Properties of ARCH Processes}
\newtheorem{defn}{Defination}
\begin{defn}
An autoregressive conditional heteroscedastic (ARCH) model with oreder p ($\ge1$) is defined as
\begin{equation} \label{e}
X_t = \sigma_t\epsilon_t\hspace{.1in}and\hspace{.1in}\sigma^2_t = c_0 + b_1X^2_{t-1} +\ldots+b_pX^2_{t-p}
\end{equation}
where $c_0\ge0$, $b_j\ge0$ are constants, $\epsilon_t \sim iid(0,1)$, and $\epsilon_t$ is independent of 
\{$X_{t-k}, k\ge1$\} for all t.
\end{defn}
The necessary and sufficient condition for (\ref{e}) to define a unique stationary process \{$X_t$\} with \\ $EX^2_t <\infty$ is 
\begin{equation} 
\sum_{i=1}^{p}b_i < 1 \end{equation}
Furthermore, for such a stationary solution, $EX_t = 0$ and 
$V(X_t) = c_0/(1-\sum_{i=1}^{p}b_i)$. 

\subsection{Estimation}
We always assume that \{$X_t$\} is a strictly stationary solution of the 
ARCH model (\ref{e}). Based on observations $X_1, X_2,\ldots, X_n$, we discuss
various methods for estimating parameters in the model. Listed below are
four types of estimators for parameters $c_0$ and $b_i$. They are the 
\textit{Conditional Maximum Likelihood Estimator}, and three \textit{Least 
Absolute Deviations Estimators}.\\ \\
\textbf{(a) Conditional Maximum Likelihood Estimator} If $\epsilon_t$ is normal in model (\ref{e}), the
negative logarithm of the (conditional) likelihood function based on observations $X_1, X_2,\ldots, X_n$,
ignoring constants, is
\begin{equation} \label{mle}
\sum_{t=p+1}^{n}(\log\sigma^2_t + X^2_t/\sigma^2_t)
\end{equation}
The (Gaussian) maximum likelihood estimators are defined as the minimizers of the function above. 
Note that this likelihood function is based on the conditional probability density function of 
$X_{p+1},\ldots, X_n$, given $X_1,\ldots,X_p$, since the unconditional probability density function,
which involves the joint density of $X_1,\ldots,X_p$ is unattainable.\\ \\
\textbf{(b) Least Absolute Deviations Estimators} The estimator discussed in (a) is derived from maximizing an approximate Gaussian likelihood. In this sense, it is an $L_2$-estimator It is well known that $L_1$-estimators are more robust with respect to heavy-tailed distributions than $L_2$-estimators. This motivates the study of various \textit{least absolute deviations estimators} for $c_0$ and $b_i$ in model (\ref{e}).\\
Now we reparametrize the model (\ref{e}) in such a way that the median of $\epsilon^2_t$, instead of the variance of $\epsilon_t$, is equal to 1 while $E\epsilon_t = 0$ remains unchanged. Under this new reparametrization, the parameters $c_0$ and $b_i$ differ from those in the old setting by a common positive constant factor. Write
\begin{equation}\label{r1}
\frac{X^2_t}{\sigma_t(\theta)^2} = 1 + e_{t1}
\end{equation}
where $e_{t1} = (\epsilon^2_t - 1)$ has median 0. This leads to the first absolute deviations estimator

\begin{equation}\label{e1}
\hat\theta_1 = \arg \min_{\theta}\sum_{t=p+1}^{n}|{X^2_t/\sigma_t(\theta)^2 -1}|
\end{equation}

\noindent which is an $L_1$ estimator based on the regression relationship (\ref{r1}).\\ \\
Alternatively, we can define another form of least absolute estimator as
\begin{equation}\label{e2}
\hat\theta_2 = \arg \min_{\theta}\sum_{t=p+1}^{n}|\log(X^2_t - \log(\sigma_t(\theta)^2)|
\end{equation}
which is motivated by the regression relationship
\begin{equation}\label{r2}
\log(X_t^2) = \log(\sigma_t(\theta)^2) + e_{t2}
\end{equation}
where $e_{t2} = \log(\epsilon^2_t)$. Hence median of $e_{t2}$ is equal to $\log\{median(\epsilon^2_t)\}$, which is 0 under the reparameterisation.\\ \\
The third L-1 estimator is motivated by the regression equation
\begin{equation}\label{r3}
X^2_t = \sigma^2_t + e_{t3}
\end{equation}
where $e_{t3} = \sigma^2_t(\epsilon^2_t - 1)$. Again under the new parameterisation, the median of $e_{t3}$ is 0. This leads to the estimator
\begin{equation}\label{e3}
\hat\theta_3 = \arg \min_{\theta}\sum_{t=p+1}^{n}|X^2_t - \sigma_t(\theta)^2|
\end{equation}
Intuitively we prefer the estimator $\hat\theta_2$ to $\hat\theta_3$ since the error terms $e_{t2}$ in regression model (\ref{r2}) are independent and identically distributed while the errors $e_{t3}$ in model (\ref{r3}) are not independent. Another intuitive justification for using $\hat\theta_2$ is that, the distribution of $X^2_t$ is confined to the nonnegative half axis and is typically skewed. Hence the log-transformation will make the distribution less skewed.\\
The minimization in (\ref{e1}) , (\ref{e2}) and (\ref{e3}) is taken over all $c_0 > 0$ and all nonnegative $b_i$'s.

\subsection{Assymptotic Properties } In this section we discuss the assymptotic properties of the estimators listed above.\\
The conditional maximum likelihood estimation remains as one of the most frequently-used methods in fitting ARCH models. To establish the assymptotic normality of the likelihood estimator some regularity conditions are required.
Let \{$X_t$\} be the unique strictly stationary solution from ARCH(p) model (\ref{e}) in which $\epsilon_t$ may not be normal. We assume that $p \ge 1$, $c_0 > 0$ and $b_i >0$ for $i = 1,2,\ldots,p$. Let $(\hat c_0, \mathbf{\hat a^T})^T$
be the estimator derived from minimizing (\ref{mle}), which should be viewed as a (conditional) quasimaximum likelihood estimator.\\
Let $\mathbf{\theta} = (c_0, \mathbf{a^T})^T$, $\mathbf{\hat\theta} = (\hat c_0, \mathbf{\hat a^T})^T$, and 
$\mathbf{U_t} = \frac{\mathit{d}\sigma^2_t}{\mathit{d}\mathbf{\theta}}$. It may be shown that $\mathbf{U_t}/\sigma^4_t$ has all its moments finite. We assume that the matrix
\[\mathbf{M} \equiv E(\mathbf{U_tU^T_t}/\sigma^4_t)\]
is positive definite. Further we assume that the errors are not very heavy tailed, ie $E(\epsilon^4_t)<\infty$. Then under the above regularity conditions, it can be established that (see Hall and Yao 2003)
\[
\frac{\sqrt{n}}{(E(\epsilon^4_t) - 1)^{1/2}}\mathbf{(\hat\theta - \theta)}\stackrel{d}\longrightarrow N(0, M^{-1})
\]
If $E(\epsilon^4_t) = \infty$ the convergence rate of $\sqrt{n}$ is no longer observable. Then the convergence rate of the likelihood estimator is dictated by the distribution tails of $\epsilon^2_t$; the heavier the tails, the slower the convergence. Moreover, the assymptotic normality of the estimator is only possible if 
$E(|\epsilon_t|^{4 - \delta}) < \infty$ for any $\delta > 0$.\\ \\
The asymptotic normality of the least absolute deviations estimator $\hat\theta_2$ in (\ref{e2}) can be established under milder conditions. To do so we will use the reparameterized model.
Let $\mathbf{\theta} = (c_0, \mathbf{a^T})^T$ be the true value under which the median of $\epsilon^2_t$ equals 1, or equivalently the median of $\log(\epsilon^2_t)$ equals 0. Define $\mathbf{U_t}$ and $\mathbf{M}$ as before. Again we assume there exists a unique strictly stationary solution \{$X_t$\} of model (\ref{e}) with $E_{\theta}(X^2_t) < \infty$. The parameters $c_0$ and $b_i$, $i = 1,2,\ldots,p$ are positive. $\mathbf{M}$ is positive definite. $\log(\epsilon^2_t)$ has median zero, and its density function \textit{f} is continous at at zero.\\
Under the above conditions, there exists a sequence of local minimizers $\hat\theta_2$ of (\ref{e2}) for which 
\[
\sqrt{n}\mathbf{(\hat\theta_2 - \theta)}\stackrel{d}\longrightarrow N(0, M^{-1}/\{4f(0)^2\})
\]
(see Peng and Yao 2003).Thus the least absolute deviations estimator $\hat\theta_2$ is asymptotically normal with convergence rate $\sqrt{n}$ under very mild conditions. In particular, the tail-weight of the distribution of $\epsilon_t$ is irrelevant as no condition is imposed on the moments of $\epsilon_t$ beyond $E(\epsilon^2_t) < \infty$\\ \\
Similar to the above result,$\sqrt{n}\mathbf{(\hat\theta_1 - \theta)}$ is also asymptotically normal with mean
\[
E[\epsilon^2_tI(\epsilon^2_t > 1) - \epsilon^2_tI(\epsilon^2_t <1)]\;[E|m_{11}|,\ldots,E|m_{(p+1)(p+1)}|]^T \;\;where\; M = (m_{ij})_{i,j}
\]
(see Peng and Yao 2003) which is unlikely to be 0. This shows that $\hat\theta_1$ is often a biased estimator.\\ \\
It can also be shown that $\sqrt{n}(\hat\theta_3 - \theta)$ is also assymptotically normal under the additional condition $EX^4_t < \infty$.
\subsection{Bootstrap in ARCH models} As indicated in the earlier section, the range of possible limit distributions for a (conditional) Gaussian maximum likelihood estimator is extraordinaily vast. In particular the limit laws depend intimately on the error distribution. This makes it impossible in heavy tailed cases to perform statistical tests or estimation based on asymptotic distributions in any conventional sense. Bootstrap methods seem the best option for tackling these problems.\\ \\
\textbf{Residual Bootstrap(m-out-of-n) for likelihood estimator:} Let \linebreak
$\tilde\epsilon_t = X_t/\sigma_t(\hat\theta)$ for $t = p+1,\ldots,n$ and let \{$\hat\epsilon_t$\} be the standardized version of \{$\tilde\epsilon_t$\} such that the sample mean is zero and the sample variance is 1. We define
\[
\hat\tau^2 = \frac{1}{n}\sum_{t=1}^{n}\tilde\epsilon^4_t - (\frac{1}{n}\sum_{t=1}^{n}\tilde\epsilon^2_t)^2
\]
Now we draw \{$\epsilon^{*}_t$\} with replacement from \{$\hat\epsilon_t$\} and
define $X^{*}_t = \sigma^{*}_t\epsilon^{*}_t$ for \linebreak $t = p+1,\ldots,m$
\;with 
\[
(\sigma^{*}_t)^2 = \hat c_0 + \sum_{i=1}^{p}\hat b_i(X^{*}_{t-i})^2
\]
and form the statistic ($\hat\theta^{*}, \hat\tau^{*}$) based on \{$X^{*}_{p+1},\ldots,X^{*}_m$\} in the same way as ($\hat\theta, \hat\tau$) based on \{$X_{p+1},\ldots,X_n$\}. It has been proved that (Hall and Yao (2003)) as $n\rightarrow \infty$, $m\rightarrow \infty$, and $m/n\rightarrow 0$, it holds for any convex set \textit{C} that
\[
\vline P\left\{ \sqrt{m}\frac{(\hat\theta^{*} - \hat\theta)}{\hat\tau^{*}}\in C|X_1,\ldots,X_n\right\} - P\left\{\sqrt{n}\frac{(\hat\theta - \theta)}{\hat\tau}\in C\right\}\vline \longrightarrow0
\]
\textbf{Weighted Bootstrap for likelihood estimator} For every $n\ge 1$, let \{$w_{nt}$\}, $t=1,\ldots,n$, be real valued row-wise exchangeable random variables independent of \{$X_t$\}. Then we define the weighted bootstrap estimators, $\hat\theta^{*}$ of $\hat\theta$ as the minimizers of
\begin{equation} \label{wml}
\sum_{t=p+1}^{n}w_{nt}[\log\sigma^2_t(\hat\theta) + X^2_t/\sigma^2_t(\hat\theta)]
\end{equation}
Under suitable regularity conditions on the weights, we can expect the consistency 
of $\hat\theta^{*}$.\\ \\
It is well known that in the settings where the limiting distribution of a statistic is not normal, 
standard bootstrap methods are generally not consistent when used to approximate the distribution
of the statistic. In particular when the the distribution of $\epsilon_t$ is very heavy-tailed in
 the sense that $E(|\epsilon_t|^d) = \infty$ for some $2< d < 4$, the Gaussian likelihood estimator is
  no longer assymptotically normal. However the least absolute deviations estimator $\hat\theta_2$ is
   assymptotically normal under very mild conditions. Hence we expect the Bootstrap methods to work
    under larger range of possible distributions for $\hat\theta_2$.\\ \\
\textbf{Weighted Bootstrap for} $\mathbf{\hat\theta_2}$ As in (\ref{wml}) we define the weighted bootstrap estimators, $\hat\theta^{*}_2$ of $\hat\theta_2$ as the minimizers of
\begin{equation}\label{wb}
\sum_{t=p+1}^{n}w_{nt}|\log(X^2_t - log(\sigma_t(\theta)^2)|
\end{equation}
Let $\sigma^2_n = V_Bw_{ni}$, $W_{ni} = \sigma^{-1}_n(w_{ni} - 1)$, where $P_B$, $E_B$ and $V_B$, respectively, denote probabilities, expectations and variances with respect to the distribution of the weights, conditional on the given data \{$X_1,\ldots,X_n$\}. The following conditions on the weights are assumed:
\begin{eqnarray}
& E_B(w_1) = 1 \label{c1} \\
& 0 < k < \sigma^2_n = o(n) \\
& c_{1n} = Cov(w_i,w_j) = O(n^{-1}) \label{c3}
\end{eqnarray}
Also assume that $\sigma^2_n/n$ decreases to 0 as $n\rightarrow\infty$. Further assume that the conditions of Result \ref{pw} hold with $U_{nj} = W_{nj}$. Then it is plausible that
\[
\vline P\{\sqrt{n}\sigma^{-1}_n(\hat\theta^{*}_2 - \hat\theta_2) \le x|X_1,\ldots,X_n\} -
P\{\sqrt{n}(\hat\theta_2 - \theta) \le x\}\vline \stackrel{P}{\longrightarrow}0 \hspace{.1in}\forall x\hspace{.01in} \in \mathbb{R}
\]

\subsection{Numerical Properties} In this section, we compare numerically the three least absolute deviation estimators with the conditional Gaussian maximum likelihood estimator for \linebreak ARCH(1) model.
 Then we check the consistency of their Bootstrap analogues.\\
We took the errors $\epsilon_t$ to have either a standard normal distribution or a standardised 
Student's $t$-distribution with $d = 3$ or $d = 4$ deegrees of freedom. We standardized the 
$t$-distributions to ensure that their first two moments are, respectively, 0 and 1. We 
took $c_0 = 1$ and $c_1 = 0.5$ in the models. Setting the sample size $n = 100$, we drew 200 
samples for each setting. We used different algorithms to find estimates for different estimation
procedures. Since the values of the parameters $c_0$ and $c_1$ estimated by the least absolute 
deviations methods differ from the numerical values specified above by a common factor
(namely the median of the square of the distribution of $\epsilon_t$), for a given sample, we 
define the absolute error as $\vline\frac{\hat c_0}{\hat c_1} - \frac{c_0}{c_1}\vline$ 
where $\hat c_0$ and $\hat c_1$ are the respective sample estimates. We average the error over all
our samples to obtain the sample average absolute error for an estimation procedure.\\
The table below displays the average absolute error for the different estimation procedures.
The first column indicates distribution of $\epsilon_t$, the second column are the estimation 
procedures, and in the third column are the corresponding average error values.
\begin{center}
\begin{tabular}{|l|||c||c|}
\hline
Distn. & Estimate & Average error\\
\hline
Normal & $\hat\theta_{ml}$ &  2.548\\
Normal & $\hat\theta_1$    &  6.936\\
Normal & $\hat\theta_2$    &  5.274\\
Normal & $\hat\theta_3$    & 16.559\\
t-3    & $\hat\theta_{ml}$ & 11.097\\
t-3    & $\hat\theta_1$    &  5.750\\
t-3    & $\hat\theta_2$    &  2.307\\
t-3    & $\hat\theta_3$    & 56.259\\
t-4    & $\hat\theta_{ml}$ & 13.107\\
t-4    & $\hat\theta_1$    &  7.054\\
t-4    & $\hat\theta_2$    &  4.528\\
t-4    & $\hat\theta_3$    & 24.253\\
\hline
\end{tabular}
\end{center}
Figures 3a), 3b) and 3c) present the boxplots for the absolute errors with error distributions being
normal, $t_3$ and $t_4$ respectively. For models with heavy-tailed 
errors, eg $\epsilon_t \sim t_d$ with $d = 3,4$ the least absolute deviation estimator 
$\hat\theta_2$ performed best. Furthermore, the gain was more pronounced when the tails were very heavy, eg $\epsilon_t \sim t_3$.From the boxplot, it can be seen that, when $\epsilon_t$ $\sim$ $t_4$, except for a few outliers, the Gaussian maximum likelihood estimator $\hat\theta_{ml}$ was almost as good as $\hat\theta_1$ and $\hat\theta_2$. However, when $\epsilon_t \sim t_3$, $\hat\theta_{ml}$ was no longer desirable. On the other hand, when the error $\epsilon_t$ was normal, $\hat\theta_{ml}$ was of course the best. In fact the absolute error of $\hat\theta_{ml}$ was larger when the tail of the error distribution was heavier, which reflects the fact that, heavier the tails are, slower is the convergence rate; see Hall and Yao (2003). However this is not the case for the least absolute deviations estimators as they are more robust against heavy tails.\\
Overall the numerical results suggest that we should use the least absolute deviations estimator 
$\hat\theta_2$ when $\epsilon_t$ has heavy and especially very heavy tails, 
eg E($|\epsilon_t|^3$) = $\infty$, while in general the Gausian maximum likelihood estimator
$\hat\theta_{ml}$ is desirable as long as $\epsilon_t$ is not very heavy-tailed.\\ \\

Next we check the consistency of the bootstrap estimators, $\hat\theta^*_{mle}$ and 
$\hat\theta^*_2$ of $\hat\theta_{mle}$ and $\hat\theta_2$ respectively. We fixed a sample of 
size 100 from the ARCH(1) process with standard normal errors, and used 200 simulations for four 
different resampling techniques: the RB ,the m-out-of-n RB and the WB. For the m-out-of-n RB, we took
$m$ to be 50. Comparing the values of $V_n$ and $V^*_n$, the results of the KS test are:
\begin{center}
{\textbf{Two-Sample Kolmogorov-Smirnov Test}}\\
Data: $V_n$ and $V^*_n$\\
Alternative hypothesis:\\
cdf of $V_n$ does not equal the cdf of $V^*_n$ for at least one sample point\vspace{0.05in}

\begin{center}
\begin{tabular}{|c||c||c||c|}
\hline
Estimate & BS Technique & KS value & p-value \\
\hline
$\hat c_{0 ml}$ & WB & 0.095 & 0.286\\
$\hat c_{1 ml}$ & WB & 0.110 & 0.152\\
$\hat c_{0 ml}$ & RB & 0.170 & 0.005\\
$\hat c_{1 ml}$ & RB & 0.125 & 0.073\\
$\hat c_{0 ml}$ & RB(m/n) & 0.1 & 0.234\\
$\hat c_{1 ml}$ & RB(m/n) & 0.095 & 0.286\\
$\hat c_{0 2}$ & WB & 0.095 & 0.286\\
$\hat c_{1 2}$ & WB & 0.130 & 0.057\\
\hline
\end{tabular}
\end{center}
\end{center}
In the table above, $\hat c_{0 ml}$ and $\hat c_{1 ml}$ denote the estimates of $c_0$ and $c_1$ respectively using the maximum likelihood estimation procedure, while $\hat c_{0 2}$ and $\hat c_{1 2}$ denote the corresponding estimates using the least absolute deviations estimator.
From the table, it can be seen that the full sample (i.e. n-out-of-n)bootstrap fails,
while m-out-of-n RB fares better. The reason that the full-sample RB fails to be consistent is 
that it does not accurately model relationships among extreme order statistics in the sample; 
see Fan and Yao 2003. WB does reasonably well for both maximum likelihood and least absolute
deviations estimation procedures.\\

\pagebreak
\begin{figure}
\begin{center}
\scalebox{.7}{\includegraphics{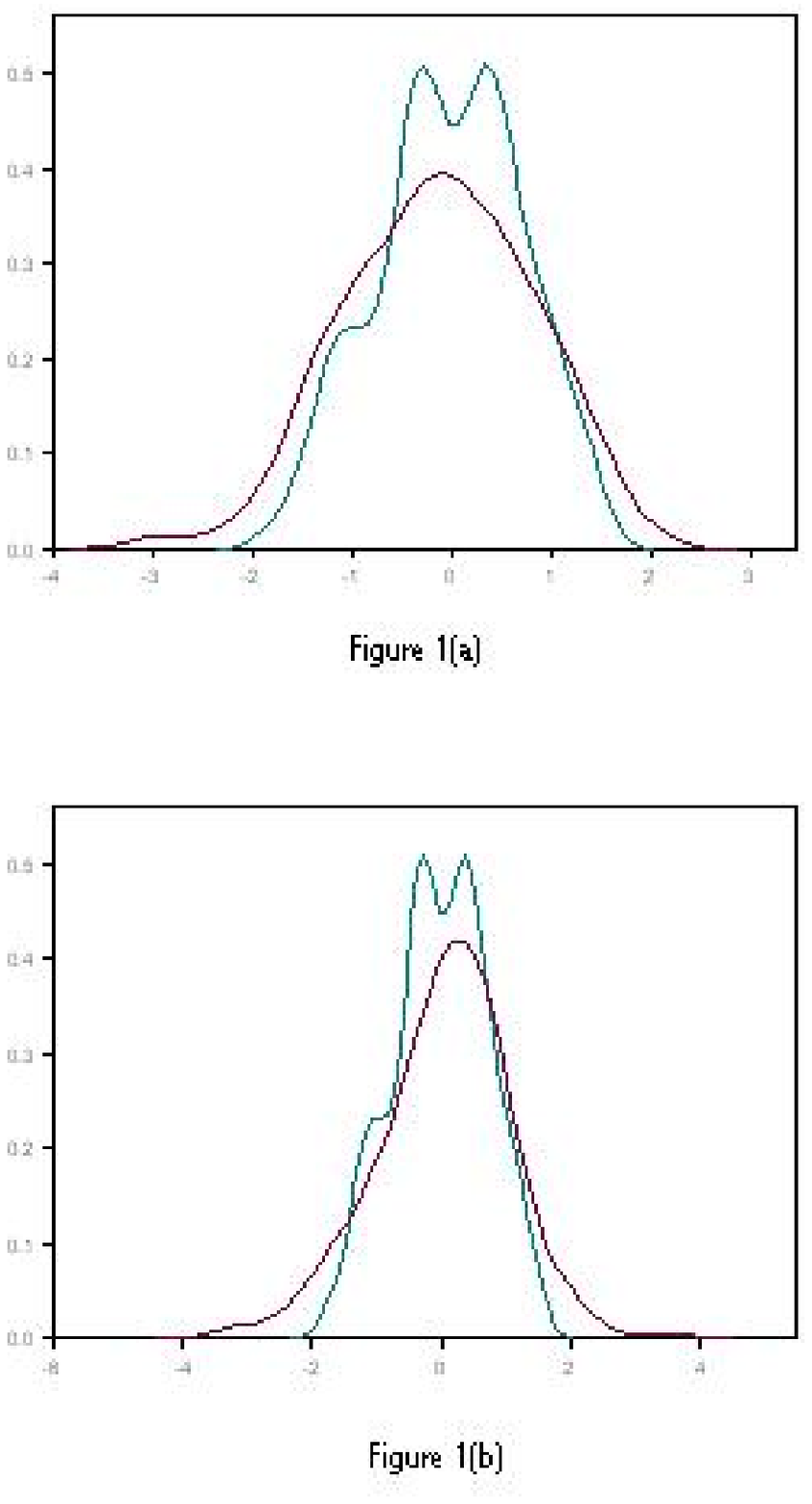}}
\end{center}
\end{figure}
\noindent\footnotesize\textbf{Figure1}: Sample density plots of $V_n$ and $V^{*}_n$ 
with $\sigma^2_1 = 1$ and $\sigma^2_2 = 2$. The green line denotes density of $V_n$, the red line for
density of $V^{*}_n$. (a) $\hat\theta^{*}_n$ is the residual bootstrap estimator, (b) $\hat\theta^{*}_n$
is the weighted bootstrap estimator.

\pagebreak
\begin{figure}
\begin{center}
\scalebox{.7}{\includegraphics{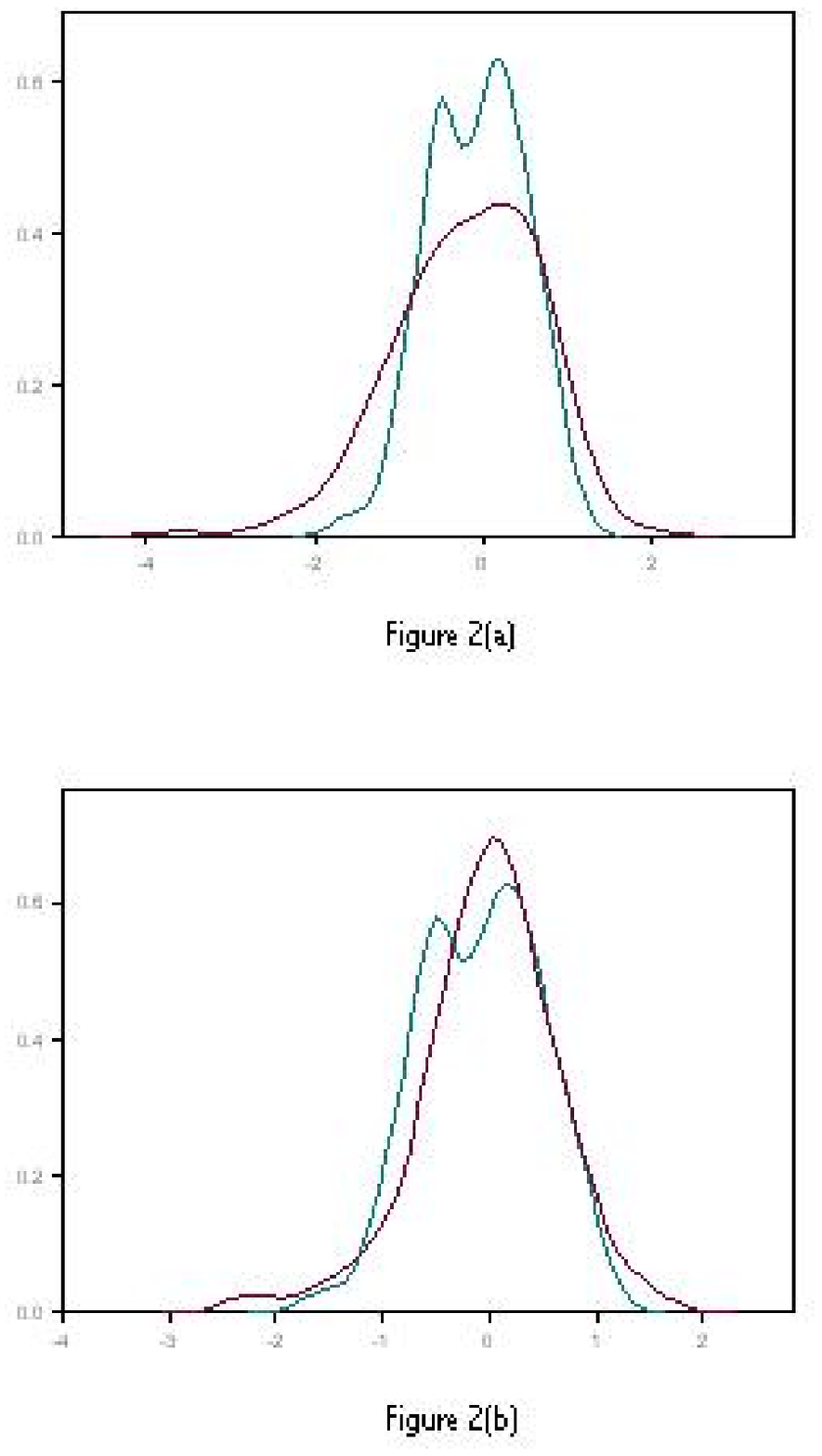}}
\end{center}
\end{figure}
\noindent\footnotesize\textbf{Figure2}: Sample density plots of $V_n$ and $V^{*}_n$ with $\sigma^2_1 = 1$ and 
$\sigma^2_2 = 10$. The green line denotes density of $V_n$, the red line for density of $V^{*}_n$. 
(a) $\hat\theta^{*}_n$ is the residual bootstrap estimator, (b) $\hat\theta^{*}_n$ is the weighted 
bootstrap estimator.

\pagebreak
\noindent\footnotesize\textbf{Figure 3}: Box plots of the absolute errors of the maximum likelihood estimates (MLE),
and the three least absolute deviations estimates (LADE). Labels 1, 2, 3 and 4 denote respectively 
the MLE, LADE1 - $\hat\theta_1$, LADE2 - $\hat\theta_2$ and LADE3 - $\hat\theta_3$. 
(a) Error $\epsilon_t$ has normal distribution, (b) Error $\epsilon_t$ has $t_3$ distribution, (c) Error $\epsilon_t$ has $t_4$ distribution.
\begin{figure}
\begin{center}
\scalebox{.6}{\includegraphics{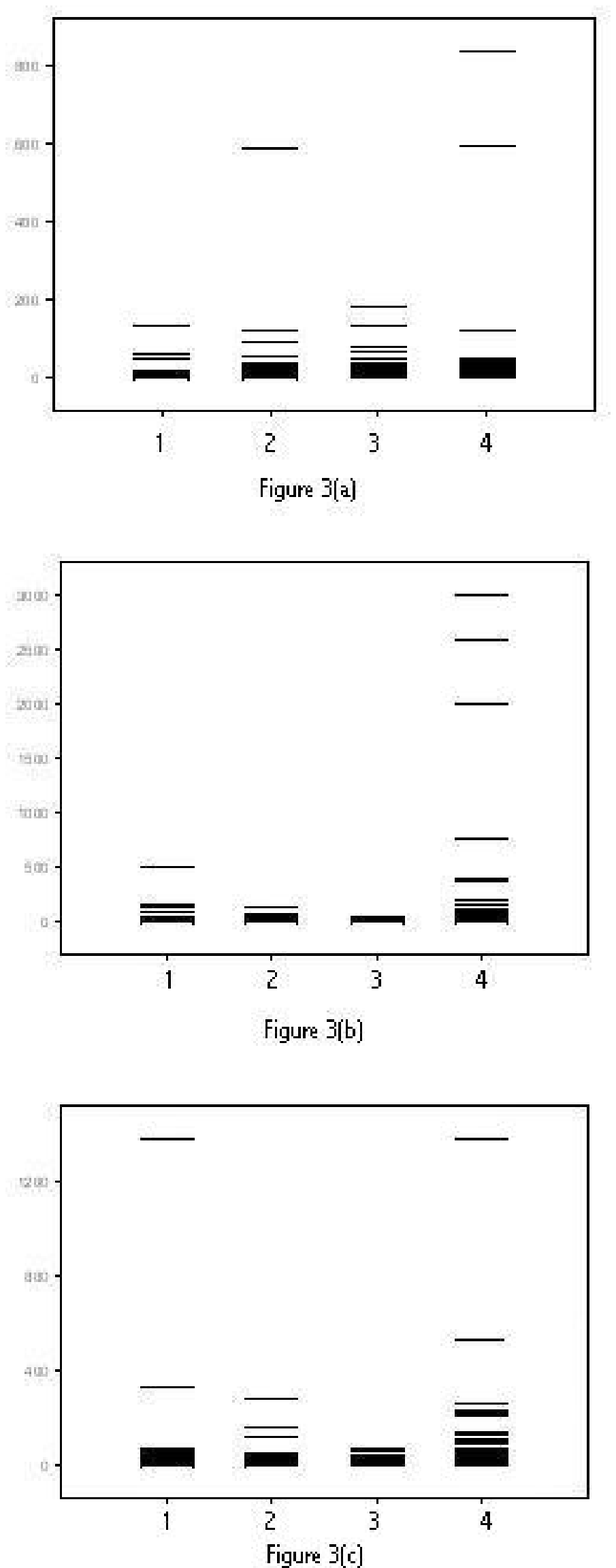}}
\end{center}
\end{figure}

\end{document}